\documentclass[11pt,letterpaper]{amsart}
\usepackage[T1]{fontenc}
\usepackage{amssymb,mathrsfs}
\usepackage{color,umoline}
\usepackage[dvipsnames]{xcolor}
\usepackage{graphicx}
\usepackage{enumitem}
\usepackage[font=small,labelfont=bf]{caption}
\usepackage{kotex}
\usepackage{tikz, caption}
\usepackage{relsize}
\usepackage[mathscr]{eucal}
\usepackage{dsfont}
\usepackage{verbatim}
\usepackage[draft]{todonotes}
\usepackage{hyperref}
\usepackage[dvipsnames]{xcolor}
\usetikzlibrary{arrows.meta}
\usepackage{subfig}

\usepackage[symbol]{footmisc}

\def\wt#1{\tilde{#1}}

\newcommand{\N}{\mathbb{N}}
\newcommand{\R}{\mathbb{R}}
\newcommand{\Z}{\mathbb{Z}}

\newcommand{\supp}{\operatorname{supp}}
\newcommand{\dist}{\operatorname{dist}}

\newtheorem{thm}{Theorem}[section]

\newtheorem{prop}[thm]{Proposition}
\newtheorem{cor}[thm]{Corollary}
\newtheorem{lem}[thm]{Lemma}
\numberwithin{equation}{section}
\newtheorem{rem}{Remark}

\newcommand\CI{{\mathcal I}}

\newcommand\CE{{\mathcal E}}
\newcommand{\cB}{{\mathcal B}}

\newcommand{\inp}[2]{\langle #1, #2\rangle}

\newcommand{\mmf}{(\mu{\tilde\mu})^\frac14}

\newcommand{\mmt}{({\mu{\tilde\mu}})^{1/2}}
\newcommand{\mm}{{\mu{\tilde\mu}}}

\newcommand{\D}{\mathcal D}

\newcommand{\sxy}{S_c(x,y)}
\newcommand{\sxyl}{2^{-l}s+\sxy}

\newcommand{\vepc}{\varepsilon_{\!\circ}}

\newcommand{\ol}[1]{\bar{#1}}

\newcommand{\smu}{\sqrt\mu}

\newcommand{\mom}{({{\tilde\mu}}/\mu)}

\newcommand{\fP}{\mathfrak P}

\newcommand{\fc}{\mathbf c}
\renewcommand{\sxy}{S_c}
\renewcommand{\sxyl}{S_c^l}

\newcommand{\cP}{\mathcal P}

\newcommand{\Be}{\begin{equation}}
\newcommand{\Ee}{\end{equation}}

\newcommand{\mU}{\mathbf U}
\newcommand{\cQ}{\mathcal Q}
\newcommand{\cD}{\mathcal D}

\newcommand{\ddk}{\partial_s^k \cP}

\newcommand{\bpqp}{{\mathrm B}_{p}(\mu, \tilde\mu)}

\newcommand{\dpp}{\delta(p,p')}
\newcommand{\sxyls}{\sxyl(x,y,s)} 

\newcommand{\eval}{2\mathbb N_0+d}

\newcommand{\su}{\sin{S_c}}
\newcommand{\cu}{\cos S_c}
\newcommand{\sul}{\sin\sxyl}
\newcommand{\cul}{\cos\sxyl}
\newcommand{\sut}{\sin^2 S_c}
\newcommand{\sult}{\sin^2\sxyl}
\newcommand{\ba}{\mathfrak A}
\newcommand{\bc}{\mathfrak G}

\newcommand{\chpl}{\chi^\pm_{\lambda, \mu}}

\newcommand{\chplp}{\chi^+_{\lambda, \mu}}

\newcommand{\chpln}{\chi^-_{\lambda, \mu}}

\newcommand{\diam}{{\rm diam }}

\newcommand{\mmfp}{\chi_{\lambda,\mu}^e\Pi_{\lambda'}\chi_{\lambda,{\tilde\mu}}^e}

\newcommand{\tmu}{\tilde\mu}
\newcommand{\ssum}[1]{\sum{}_{{}_{\mathlarger #1}}}

\newcommand{\nb}{\mathcal N_{\bf e}}

\newcommand{\nc}{\mathcal N_{\bf c}}

\newcommand{\wmu}{\tilde{\mu}}
\newcommand{\dO}{\mathcal O_{\!\ast}}

\begin{document}

\author[E. Jeong]{Eunhee Jeong}
\address[Jeong]{Department of Mathematics Education and  Institute of Pure and Applied Mathematics, Jeonbuk National University, Jeonju 54896, Republic of Korea}
\email{eunhee@jbnu.ac.kr}

\author[S. Lee]{Sanghyuk Lee}
\address[Lee]{Department of Mathematical Sciences and RIM, Seoul National University, Seoul 08826, Republic of  Korea}
\email{shklee@snu.ac.kr}

\author[J. Ryu]{Jaehyeon Ryu}
\address[Ryu]{School of Mathematics, Korea Institute for Advanced Study, Seoul 02455, Republic
of Korea}
\email{jhryu@kias.re.kr}

\title[Hermite operator]
{
Endpoint eigenfunction bounds for 
\\
 the Hermite operator
}

\begin{abstract}  
We establish  the optimal $L^p$, $p=2(d+3)/(d+1),$ eigenfunction bound for  the Hermite operator $\mathcal H=-\Delta+|x|^2$ on $\mathbb R^d$. 
Let  $\Pi_\lambda$ denote  the  projection operator to the vector space  spanned by the eigenfunctions of $\mathcal H$ with eigenvalue $\lambda$. 
The optimal $L^2$--$L^p$ bounds on $\Pi_\lambda$, $2\le p\le \infty$, have been known by the works of Karadzhov and  Koch-Tataru
except  $p=2(d+3)/(d+1)$.  For $d\ge 3$,  we prove the optimal  bound for the missing endpoint case. Our result is built  on 
a new phenomenon: improvement of  the bound due  to asymmetric localization near the sphere $\sqrt\lambda \mathbb S^{d-1}$. 
\end{abstract}

\maketitle

\section{Introduction}  
The Hermite operator $\mathcal H=-\Delta+|x|^2$ on $\mathbb R^d$ has a discrete spectrum $\lambda\in 2\mathbb N_0+d$, $\mathbb N_0:=\mathbb N\cup \{0\}$. 
For $\alpha\in \mathbb N_0^d$, we denote by  $\Phi_\alpha$ the $L^2$--normalized Hermite function,  which is  an  eigenfunction of $\mathcal H$ with eigenvalue $2|\alpha|+d$.  
The set  $\{\Phi_\alpha: \alpha\in \mathbb N_0^d\}$  forms an orthonormal basis in $L^2$.    
 Let  $\Pi_\lambda$ denote  the spectral projection operator  to the vector space  spanned by the eigenfunctions with  eigenvalue $\lambda$, i.e.,  
\begin{align*}
  \textstyle  \Pi_\lambda f =\sum_{\alpha:\, d+2|\alpha|=\lambda}\, \langle f,\Phi_\alpha \rangle \Phi_\alpha. 
\end{align*}

In this paper, we are concerned with  bounds on  the operator norm 
$ \| \Pi_\lambda  \|_{2\to q}$ for $2\le q\le \infty$, 
 where $\| T  \|_{s\to r}$ denotes the norm of an  operator $T$ from $L^s$ to $L^r$. 
The sharp bound in terms of $\lambda$ has  been of interest in connection to Bochner-Riesz summability of the Hermite expansion. 
See, Askey and Wainger \cite{AW},  Karadzhov \cite{K94}, and Thangavelu \cite{Th98} (also, see \cite{CLSY,CLWY} and \cite{LR} for recent developments).  The bounds independent of $\lambda$ have  applications to the strong unique continuation problem for  parabolic equations  \cite{e00, ev01, T09, JLR-heat}.

The bounds on  $\|\Pi_\lambda \|_{2\to q}$ have been almost completely understood.  
 When $d=1$,  the sharp bounds follow from  those  on $L^p$ norm of   the 
Hermite functions in $\mathbb R$ (\cite{Th93}).  Let 
\[q_\circ=\frac{2(d+3)}{d+1}.\]
In higher dimensions $d\ge 2$, by the works of Karadzhov \cite{K94} and Koch-Tataru \cite{T05}  it is known  that,  for $q\in [2,\infty]\setminus \{q_\circ\}$, 
\Be
\label{main}
\|\Pi_\lambda \|_{2\to q} \sim \mathrm B_q(\lambda) :=\max\big( \lambda^{-\frac12+\frac d2\delta(2,q)}, \  \lambda^{\frac d6\delta(2,q)-\frac16}, \ \lambda^{-\frac12\delta(2,q)}\big),
\Ee
where  $\delta(r,s)=r^{-1}-s^{-1}$.   
(See {\it Notation} below for the precise meaning of $\sim$.) The bound for $q=2$ is clear from Bessel's inequality, and that for $q=\infty$ is  a consequence 
of  the estimate for the kernel of $\Pi_\lambda$ (see \cite[Lemma 3.2.2]{Th93}).  Karadzhov \cite{K94} showed 
$    \|\Pi_\lambda \|_{2\to {2d}/(d-2)}\le C$ for a constant $C$.  Thangavelu  \cite{Th98} considered a local estimate  over a compact set $K$ and he obtained a sharp bound on $\|\chi_K\Pi_\lambda \|_{2\to 2(d+1)/(d-1)}$.  A systematic study  was carried out by Koch and Tataru  and they almost completely characterized $L^2$--$L^q$ bounds  (\cite[Corollary 3.2]{T05}) including the lower bounds $\|\Pi_\lambda \|_{2\to q} \ge C\,\! \mathrm B_q(\lambda)$ for some constant $C>0$ when $2<q<\infty$ (see \cite[Section 5]{T05}). 
 
  However, prior to the present work,   the  optimal estimate   remains unsettled  for  $q=q_\circ$.  
  This contrasts with the spectral projections of other related differential operators whose  optimal $L^2$--$L^q$ bounds are well understood  \cite{sogge88, SZ,  K07, JLR-twist}. 
     By virtue of  localized estimates over annuli  (see \eqref{kt-nu} below),  it was known  \cite{T07}   that
    \[     \|\Pi_\lambda \|_{2\to q_\circ }\le C\lambda^{-\frac{1}{2(d+3)}}(\log\lambda)^{1/q_\circ}. \] 
 When $d=1$,  the estimate  fails  without the logarithmic factor. However,  when $d\ge 2$,  it was conjectured  in \cite{T05} that the natural bound \eqref{main} extends to the missing endpoint  $q=q_\circ$. This  case is the most significant since   interpolation recovers the sharp bounds for $ 2<q<2d/(d-2)$. 

We prove the conjecture is true for every  $d\ge 3$. 

\begin{thm}\label{endpoint}
Let $d\ge 3$.  Then, 
\[
 \|\Pi_\lambda \|_{2\to q_\circ}\sim \lambda^{-\frac{1}{2(d+3)}}.
\]
\end{thm}

It is likely that the theorem continues to be true for $d=2$ but our argument in this paper is not enough to prove this case. 

\subsubsection*{Localized estimate} For $ \mu \in \mathbb D^{-}\!:= \{ 2^k:-k\in \mathbb N\}$, we set 
\[   
A^{\pm}_{\mu}=\{x:   \pm(1- |x|)\in [\mu, 2\mu]\}, \quad  \  A^{\pm}_{\lambda, \mu}=\big\{x:   \lambda^{-\frac12} x\in A^{\pm}_{\mu} \big\}, 
\]
respectively. 
For simplicity, we also denote 
\[\chi^\pm_\mu = \chi_{A^{\pm}_\mu}, \qquad \chi^\pm_{\lambda, \mu}=\chi_{A^{\pm}_{\lambda,\mu}}.\]
Of special interest is the estimate over the region near the sphere  $ \sqrt\lambda \mathbb S^{d-1}=\{x :|x|=\sqrt\lambda\,\} $, across which  the kernel of  $\Pi_\lambda$ exhibits different behaviors.  
Koch and Tataru \cite{T05} considered  the localized operator $\chpl\Pi_\lambda$. They proved  the following sharp bounds: 
\begin{equation}  
\label{kt-nu}
 \|\chplp\Pi_\lambda \|_{2\to q}\le  C
                                      \begin{cases}
\lambda^{-\frac12\delta(2,q)}\mu^{\frac14-\frac{d+3}{4}\delta(2,q)}, &   \  \  \, 2\le q\le \frac{2(d+1)}{d-1},
\\[2pt]
(\lambda\mu)^{\frac d2\delta(2,q)-\frac12}, &   \   \ \, \frac{2(d+1)}{d-1} \le q\le \infty,
                                      \end{cases}
\end{equation}
for $\lambda^{-\frac23}\le\mu\le 1/4$  \cite[Theorem 3]{T05}.\footnote{This is what was  proved in \cite{T05}, where some  care with the notation $\ell_\lambda^\infty L^p$ seems necessary.} 
Summation over $\mu$ and interpolation with the previously known bound   give  \eqref{main} except for  $q=q_\circ$. 
Meanwhile, the estimates for $\chpln\Pi_\lambda $ are of less interest,   since $\chpln\Pi_\lambda $  has much smaller bounds thanks to rapid decay of its kernel 
(e.g.,  see \eqref{est:ext--} below).

Let   $\mu_\circ =\lambda^{-2/3}$ and $c_\circ= (100d)^{-2}$.   Thanks to the estimates in \cite{T05} which are mentioned above, we already have the desired $L^2$--$L^{q_\circ}$ bounds on 
$  \sum_{ c_\circ\le \mu<1}\chi^{+}_{\lambda,\mu}  \Pi_\lambda $ and $( 1- \sum_{ \mu_\circ<\mu< 1}\chi^{+}_{\lambda,\mu}) \Pi_\lambda$. (We refer the reader forward to Section \ref{pfpf} for the detail.)
Therefore,  to prove  Theorem \ref{endpoint} 
 it is sufficient to consider $L^2$--$L^{q_\circ}$ bound 
 on the operator 
\[\textstyle \Pi_\lambda':= \sum_{ \mu_\circ<\mu< c_\circ}\chi^{+}_{\lambda,\mu}  \Pi_\lambda.\]  
By duality, 
$\|\Pi_\lambda'\|_{2\to q_\circ}^2=\| \Pi_\lambda'(\Pi_\lambda')^*\|_{q_\circ'\to q_\circ}$. So, if one  shows   $\| \Pi_\lambda'(\Pi_\lambda')^*\|_{q_\circ'\to q_\circ}\le C   \lambda^{-\delta(q_\circ',q_\circ)/2}$, 
\eqref{endpoint} follows since the lower bound  is already shown in \cite{T05}. 
 Since  $\Pi_\lambda^\ast=\Pi_\lambda$ and $\Pi_\lambda^2=\Pi_\lambda$,  we can write   
 \[ \textstyle  \Pi_\lambda'(\Pi_\lambda')^*  = \sum_{ \mu_\circ<\mu, \tilde \mu< c_\circ}  \chi_{\lambda,\mu}^+ \Pi_\lambda   \chi_{\lambda,{\tilde\mu}}^+. \]
  This naturally leads us to consider  $L^p$--$L^q$ bounds on the operators $\chplp \Pi_\lambda  \chi_{\lambda,{\tilde\mu}}^+$  for general exponents $p,q$, not necessarily restricted to the case $p=q'$.  
Note that  $\|\chplp \Pi_\lambda  \chi_{\lambda,{\tilde\mu}}^+\|_{p\to q}\le \|\chplp \Pi_\lambda\|_{2\to q} \|\chi_{\lambda,{\tilde\mu}}^+  \Pi_\lambda \|_{2\to p'}.$ 
The bound \eqref{kt-nu} yields 
\Be
\label{notimproved}
   \|\chplp \Pi_\lambda  \chi_{\lambda,{\tilde\mu}}^+\|_{p\to q}  \le C \lambda^{-\frac12\delta(p,q)}  (\mm)^{\frac14-\frac{d+3}8\delta(p,q)}    
   \Ee
   for $2(d+1)/(d+3)\le p\le 2\le q\le 2(d+1)/(d-1)$.    Attempting to add those  bounds with some interpolation trick  does not seem feasible  to recover the missing endpoint case.  Due to  the optimality of  \eqref{kt-nu}, the bound \eqref{notimproved} with $p=q'$ can not be improved if $\mu\sim \tilde \mu$. 
However, this does not exclude the possibility of an improved bound when $\mu\not\sim \tilde \mu$. 

  \subsubsection*{Asymmetric  localization} Our main novelty is in the following theorem which shows improvement of the bounds thanks  to the
 \emph{asymmetric} localization near $ \sqrt\lambda \mathbb S^{d-1}$. In other words, the bound on the operator $\chplp \Pi_\lambda  \chi_{\lambda,{\tilde\mu}}^+$, $\tilde\mu\le \mu$, compared with \eqref{notimproved},  significantly improves
as ${\tilde\mu}/\mu$ gets smaller.

\begin{thm}
\label{unb} Let $d\ge 3$,  $\mu,{\tilde\mu}\in \mathbb D^{-}\!$, and $\lambda^{-2/3}\le {\tilde\mu}\le\mu \le c_\circ= (100d)^{-2}$. If $2< q\le 2(d+1)/(d-1)$, then 
there are positive  constants  $c, C$, independent of $\mu, \tilde\mu$, and $\lambda$,  such that 
\Be
\label{unbalanced}    \|  \chi_{\lambda,\mu}^+\Pi_\lambda  \chi^{+}_{\lambda, {\tilde\mu}}  \|_{q'\to q} \le C 
\lambda^{-\frac12\delta(q',q)}(\mu{\tilde\mu})^{\frac14-\frac{d+3}{8}\delta(q',q)}  \mom^{c}.
\Ee
\end{thm}
 
 This is a new phenomenon which has not been observed before.  
 Our approach in this paper provides an elementary alternative proof of  the estimate \eqref{kt-nu} which corresponds to the case $\mu\sim \tilde\mu$. However, 
 as we shall see later,    to obtain the improved bound  \eqref{unbalanced}  is far less trivial (Section \ref{sec:endpoint} and \ref{sec:l2}).  The estimate \eqref{unbalanced} can also be extended  to some $p, q$ other than $p=q'$ by interpolation (see also \cite{JLR-hermite}). 
  
  We briefly explain how  one can obtain the missing endpoint bound from the estimate \eqref{unbalanced}.  
More details are to be provided in Section \ref{pfpf}.    We write 
\[
 \Pi_\lambda'(\Pi_\lambda')^*=
 \sum_{k} \mathcal T_k:=
 \sum_{k} \sum_{(\mu, {\tilde\mu})\in  \mathfrak D_k }   \chi_{\lambda,\mu}^+ \Pi_\lambda   \chi_{\lambda,{\tilde\mu}}^+,
 \]
 where $\mathfrak D_k=\{(\mu, \tilde\mu):  {{\tilde\mu}}/\mu\in [ 2^k, 2^{k+1}), \mu, \tilde\mu\in (\mu_\circ, c_\circ)\}$ and $\mu, \tilde\mu\in \mathbb D^{-}$. 
Considering the adjoint operators,  we also have  improved bounds with the additional factor $(\mu/\tilde \mu)^c$ when $\mu<\tilde\mu$.
Thus,  applying Theorem \ref{unb} with $q=q_\circ$, one  gets  $\| \mathcal T_k\|_{q_\circ'\to q_\circ} \le C   2^{-c|k|} \lambda^{-\delta(q_\circ',q_\circ)/2}$, which 
 consequently shows the desired endpoint bound $\|\Pi_\lambda' \|_{2\to q_\circ}\lesssim \lambda^{-{1}/{2(d+3)}}$.

\subsubsection*{Weighted estimates}
Through the same argument we can obtain a more general result which contains  the endpoint bound in Theorem \ref{endpoint}. In fact,  using Theorem \ref{unb}, we prove the following weighted estimates  which were conjectured in \cite[Remark 3.1]{T05}. Let  us set 
\[ w_\pm(x)=1+ \lambda^{-1/3}(\lambda-|x|^2)_\pm. \]

\begin{cor} 
\label{cor:weighted} Let $d\ge 3$ and $2< q\le \infty$.  Set $\gamma=\gamma(q):=\min\big(\frac{d+3}4\delta(2,q)-\frac14,  \frac12- \frac{d}2\delta(2,q) \big) $. 
Then, for  $N>0$, there is a constant $C=C(N)$ such that  
\Be 
\label{weighted}
  \big\|   w_+^{\gamma}  w_{-}^N\, \Pi_\lambda f \big\|_{q}\le C \lambda^{\frac d6\delta(2,q)-\frac16} \|f\|_2 .
  \Ee
\end{cor}

In particular, if we take $q=q_\circ$, $\gamma=0$ and hence $w_+^{\gamma}  w_{-}^N\ge 1$.  
So,   the endpoint  bound $\|\Pi_\lambda \|_{2\to q_\circ}\le C \lambda^{-{1}/{2(d+3)}}$  follows from the estimate  \eqref{weighted}.

\smallskip

\noindent {\it Organization.}  In Section 2,  we  obtain  some preparatory results, and we prove Corollary \ref{cor:weighted} while assuming  Theorem \ref{unb}. 
  In Section 3, we reduce the proof of Theorem \ref{unb}  to 
that of an $L^2$ estimate, which we show in Section 4. 

\smallskip

\noindent {\it Notation.} For  nonnegative quantities   $A$ and $B$,    $B \lesssim A$ means there is a constant $C$, depending only on dimensions, such that 
$B\le CA $.  Likewise,   $A\sim B$ if and only if   $B \lesssim A$ and $A \lesssim B$. By $D=O(A)$  we denote 
$|D|\lesssim A$.  Additionally, we denote $A\gg B$ if there is a sufficiently large constant $C>0$ such that $A\ge CB$.  By   $c$ and $\vepc$ we denote  positive  constants which are chosen to be small enough.

\section{The projection operator $\Pi_\lambda$}
In this section,  we make reductions toward  the proof of  Theorem \ref{unb}. We also obtain  some estimates for  the projection operator $\Pi_\lambda$, which are to be used in  Section \ref{sec:endpoint} and  \ref{sec:l2}. 
At the end of this section, we provide the proof of  Corollary \ref{cor:weighted}  while assuming    Theorem \ref{unb}. 
 
\subsection{The kernel of  $\Pi_\lambda$}  The  Hermite-Schr\"odinger propagator $e^{-it\mathcal H} f$ is given   by 
\Be 
\label{schro-op}
 \textstyle e^{-it\mathcal H} f=\sum_{\lambda\in \eval} \, e^{-it\lambda} \Pi_\lambda f, \quad  f\in \mathcal S(\mathbb R^d).
  \Ee
  Here, $\mathcal S(\R^d)$ denotes the  Schwartz class on $\R^d$.
 It is easy to  see the series converges uniformly.  
 Since $\mathcal H^N \Phi_\alpha=(2|\alpha|+d)^N \Phi_\alpha$,  integration by parts  gives $ \inp {f}{\Phi_\alpha}=(d+2|\alpha|)^{-N} \inp {\mathcal H^N f}{\Phi_\alpha}.$ Note $\|\Phi_\alpha\|_p\le C(1+|\alpha|)^{{d}/4}$  for $1\le p\le \infty$ (e.g., see  \cite{Th93}).  
 Thus, $ | \inp {f}{\Phi_\alpha}|\le C (d+2|\alpha|)^{-N+d/4}
 $  since $\mathcal H^N f\in \mathcal S(\mathbb R^d)$. Therefore, taking $N$ large enough, we see the series converges uniformly.

Note  $\frac1{2\pi}\int_{-\pi}^\pi e^{i\frac t2(\lambda-\lambda')}  dt =\delta(\lambda-\lambda'),$ $\lambda, \lambda'\in  2\mathbb N_0+d$.  So, we have
\[  \Pi_\lambda f =  \frac1{2\pi} \int_{-\pi}^\pi    \sum_{\lambda'\in\eval}  
e^{i\frac t2(\lambda-\lambda')} \Pi_{\lambda'}\! f\, dt, \quad  f\in \mathcal S(\mathbb R^d),\] 
 since the series converges uniformly. By \eqref{schro-op} it follows that
  \[
 \Pi_\lambda f = \frac{1}{2\pi}\int_{-\pi}^{\pi} e^{i\frac t2(\lambda-\mathcal H)}fdt, \quad   f\in \mathcal S(\mathbb R^d).
\]
Now, combining this and   Mehler's formula which expresses the kernel  of $e^{-it\mathcal H}$ 
(e.g.,  see \cite[p.11]{Th87}   and  \cite{SoTo}), we get the following.

\begin{lem}
\label{kernel}
Let $\lambda\in 2\mathbb N_0+d$. Set $ \mathfrak a(t) =(2\pi i\:\sin t)^{-d/2}{e^{i \pi d/4}}$ and 
\begin{align*}
\phi_\lambda (x,y,t)=\frac{\lambda t}{2}+\frac{|x|^2+|y|^2}{2}\cot t-  \inp xy \csc t\,. 
\end{align*}Then for $f\in \mathcal S(\mathbb R^d)$,  we have 
\Be
\label{proj-int}
    \Pi_\lambda f (x)=   \frac{1}{2\pi}  \int_{-\pi}^{\pi}   \mathfrak a(t) \int   e^{i\phi_\lambda(x,y,t)}   f(y)\, dy\,dt. 
   \Ee
\end{lem}

In what follows,  by $T(x,y)$ we denote  the kernel of an operator $T$.
Note  
\Be 
\label{rot}\phi_\lambda(\mU x,\mU y,t)=\phi_\lambda(x,y,t), \quad  \mathbf U\in \mathrm O(d),
\Ee 
 where $\mathrm O(d)$ denotes the orthonormal group in $\mathbb R^d$. Obviously, 
\eqref{rot} implies  $ \Pi_\lambda (\mU x, \mU y)=\Pi_\lambda(x,y) . $
Kochneff \cite{kochneff} made the same observation by using  the properties of the Hermite functions.

\medskip

\noindent\emph{Dyadic decomposition.} 
We dyadically  decompose  the integral  \eqref{proj-int} away from the singularities $0, \pm\pi$ of  $\mathfrak a$. 
To do so,  let $\psi\in C^\infty_c([\frac14,1])$ such that $\sum_{j\in\Z} \psi(2^j t) =1$ for $t>0$, and then  define 
$ \psi^0$ by
\begin{equation}
\label{j-decomp}
 \psi^0(t)+ \sum_{j\ge 4} \big( \psi(2^j t)+ \psi(-2^j t)  +\psi(2^j(t+\pi))+\psi(2^j(\pi-t))\big)=1
\end{equation}
 for $t\in(-\pi,\pi)\setminus\{0\}$. So, $ \supp\psi^0\subset (-\pi, \pi)\setminus\{0\}$.
   For a bounded  function $\eta$ supported in $[-\pi, \pi]$, we consider  
   \Be \label{proj-pi}
\Pi_\lambda[\eta]= \int \eta(t)e^{i\frac t2(\lambda-\mathcal H)}dt.
\Ee
  Since  $\|e^{-it\mathcal H} f\|_2=\|f\|_2$ for $t\in \mathbb R$, we have   
  \Be
  \label{easy-l2} 
  \|\Pi_\lambda[\eta]\|_{2\to 2} \le \|\eta\|_1.
  \Ee
  
For simplicity,   let us denote 
\[ \psi_j^\pm(t) = \psi(\pm 2^j t), \qquad \psi_j^{\pm \pi} (t)=\psi(2^j(\pi\mp t)), 
\qquad 
\psi_j^0= 
\begin{cases} \psi^0,  &j=4 
\\[0pt]
\ 0,  &{ j>4}
\end{cases}\, .\]
Then, using \eqref{j-decomp} we  decompose 
\Be
\label{decompj}
 \Pi_\lambda f =  \sum_{\kappa=0, \pm, \pm \pi}  \, \sum_{j\ge 4}   \Pi_\lambda [\psi_j^\kappa ]f . 
\Ee
Validity of the decomposition  is clear   since the sum on the right hand side 
converges to $\Pi_\lambda$ as a  bounded operator on $L^2$.  Indeed, \eqref{easy-l2} gives
$\| \Pi_\lambda [\psi_j^\kappa ]\|_{2\to 2}\lesssim 2^{-j}$, $\kappa= \pm, \pm \pi,$  so the convergence follows.  

We obtain the estimate for $\Pi_\lambda$ by considering each $\Pi_\lambda [\psi_j^\kappa ]$. 
However, thanks to  symmetric properties of the kernels $\Pi_\lambda [\psi_j^\kappa ](x,y)$,  the matter reduces to showing 
the estimates for $ \Pi_\lambda [\psi_j^+]$ and $\Pi_\lambda[\psi^0]$. Indeed, observe  $ \phi_\lambda(x,y,-t)=-\phi_\lambda(x,y,t)$ and $\phi_\lambda(x,y,\pm (\pi-t))=\pm\lambda \frac \pi 2 \mp\phi_\lambda(x,-y,t).$
Then, changes of variables give
\begin{align}
\label{symmetric} 
\Pi_\lambda[\psi_j^{-}](x,y)&=C_d \,\overline{\Pi_{\lambda}[\psi_j^{+}]}(x,y), 
\\[4pt]
\label{symmetric2} 
\Pi_\lambda[\psi_j^{\pm \pi}](x,y)&= C_d' \, \Pi_{\lambda} [\psi_j^{\mp}](x,-y),
\end{align}
where  $C_d$, $C_d'$ are constants satisfying $|C_d|=|C_d'|=1$. This implies 
\[ \textstyle \|\ssum{j}\, \Pi_\lambda[\psi_j^{+}]\|_{p\to q}= \|\ssum{j}\, \Pi_{\lambda} [\psi_j^{\kappa}]\|_{p\to q},  \qquad   \kappa=-,\, \pm \pi.\]  
 
\subsubsection*{Rescaled operators} 
Instead of $\Pi_\lambda$ and $\Pi_\lambda[\eta]$, it is more convenient to work with the rescaled operators  $\fP_\lambda$,  $\fP_\lambda[\eta]$ whose kernels are  given by  
\[
\fP_\lambda(x,y)=\Pi_\lambda(\sqrt \lambda x, \sqrt \lambda y), 
\quad  \  \
 \fP_\lambda[\eta](x,y)= \Pi_\lambda[\eta](\sqrt \lambda x, \sqrt \lambda y),
\]
 respectively. By rescaling we have,   for any measurable sets $E,$ $F$, 
\begin{equation}
\label{eq:norm-scaled}
  \|  \chi_{E} \mathfrak \fP_\lambda [\eta]  \chi_{F} \|_{p\to q} =  \lambda^{\frac d2(\frac1p-\frac1q-1)}\|  \chi_{\sqrt\lambda E} \Pi_\lambda [\eta] \chi_{\sqrt\lambda F} \|_{p\to q}.
   \end{equation}

 To prove Theorem \ref{unb},   we need only  to consider the case  
\Be\label{mucon}   \lambda^{-2/3}\lesssim {\tilde\mu},  \mu  \le  c_\circ,  \quad   \wmu/\mu\le  \varepsilon_\circ \Ee 
for a small  constant $\varepsilon_\circ>0$, since \eqref{unbalanced}  follows from \eqref{kt-nu} if  $\mu\sim {\tilde\mu}$.

\begin{prop} 
\label{improvedbound}  Let $d\ge 3$ and $\mu, \wmu$ satisfy \eqref{mucon}.   Set  \[ \mathrm B_{p}(\mu, {\tilde\mu}) = \lambda^{\frac{d-1}2\delta(p,p')-\frac d2}  (\mm)^{\frac14-\frac{d+3}8\delta(p,p')}    .\] 
 Suppose $0<\dpp< \min(2/(d-1), 2/3)$.  Then,  for a positive constant  $c$ independent of $\mu, \tilde \mu$,  we have  
\Be
\label{scaled-c}
 \|  \chi_\mu^+\fP_\lambda  \chi^{+}_{{\tilde\mu}}  \|_{p\to p'} \lesssim  
\bpqp  \mom^{c}.  
\Ee
\end{prop}
 
Note that $
 \min(2/(d-1), 2/3)> 2/(d+1)$ when $d\ge 3$. Thus, Theorem \ref{unb} follows. 
 In fact,  \eqref{unbalanced}  holds in a sightly bigger range. 
The rest of the paper is  devoted to the proof of Proposition \ref{improvedbound}. 

\subsection{The phase function $\mathcal P$} 

Let us set 
\begin{align}
\label{p-def} 
\cP(x,y,s)=\frac{s}{2}+\frac{|x|^2+|y|^2}{2}\cot s-  \inp xy \csc s.
\end{align}

Note  $\cP(x,y,s)=\phi_1(x,y,s)$.  By Lemma \ref{kernel} we have 
\begin{align}\label{S3-ptkernel}
\fP_\lambda[\eta] (x,y)
    =\int   (\eta\mathfrak a )(s)  e^{i\lambda \mathcal{P}(x,y,s)}ds,  
    \end{align}
which we shall extensively  make use of throughout the paper. 
To estimate the kernels,  we have a close look at  the phase function $\mathcal P$.    
A  computation together  with an elementary trigonometric identity gives 
\begin{align}
\label{ph-d}
   \partial_s \mathcal P(x,y,s) 
           =-\frac{\mathcal{Q}(x,y,\cos s)}{2\sin^2s},  
\end{align}
where 
\Be
\label{q-def}
\begin{aligned}
  \mathcal{Q}(x,y, \tau)
&:=(\tau- \inp xy)^2-\mathcal D(x,y),
\\ 
\mathcal D(x,y)&:=1+\inp xy^2-|x|^2-|y|^2.
\end{aligned}
\Ee

The stationary  points of $\mathcal P(x,y,\cdot)$ are  given by the zeros of $\cQ(x,y,\cos \cdot)$.   
So,  $\cD(x,y)$, which determines  the nature of those stationary points,  plays a significant role in estimating   the 
kernel $\fP_\lambda[\eta] (x,y)$.   Note 
\begin{align}
 \label{angle} 
\mathcal D(x,y)=-|x|^2|y|^2 \sin^2\!\!\measuredangle(x,y) +(1-|x|^2)(1-|y|^2),
 \end{align}
 where $\measuredangle(x,y)$ denotes  the angle  between $x$ and $y$. 
 Since $(1-|x|^2)(1-|y|^2)\sim \mm$ for $(x,y)\in  A^+_\mu\times {A}^+_{{\tilde\mu}}$,   
 one can control  $\cD(x,y)$ by 
the relative size of  $\measuredangle(x,y)$  to $\mmt$.

\subsection{Preliminary decomposition}  Fix a constant $c>0$ such that $1/(20d)$ $\le c\le 1/(10d)$. 
We partition the unit sphere $\mathbb S^{d-1}$ into finite disjoint subsets  $\{ S_j\}$  of diameter less than $c$.   
Then, set $A_j =\{ x\in A^+_\mu:   |x|^{-1}x\in S_j\}$ and   $\tilde A_j =\{ x\in  A^+_{{\tilde\mu}}:   |x|^{-1}x\in S_j\}$, which respectively partition  the annuli  $A^+_\mu$   and $A^{+}_{{\tilde\mu}}$ into disjoint sets
of diameter $\le c$. So,  we have 
\Be
\begin{aligned}
\label{coarse-decomp}
A_\mu^+=\bigcup{}_{{}_{\mathlarger k}}\, A_k,   \qquad
   A_{{\tilde\mu}}^{+}=\bigcup{}_{{}_{\mathlarger l}}\,  \tilde A_l.
\end{aligned}
 \Ee
This and \eqref{decompj} give  
 \begin{equation}
\label{decomp-proj2}
\chi_\mu^+ \fP_\lambda\chi_{\tilde \mu}^+ f =  \sum_{k,l} \sum_{\kappa=0, \pm, \pm \pi}   \sum_{j\ge 4}   \chi_{A_k}  \fP_\lambda [\psi_j^\kappa ] \chi_{{\tilde A_l}} f. 
 \end{equation}
 
  Using this coarse decomposition, we can distinguish minor parts whose contributions are negligible. 
   More precisely, we have the following.

\begin{lem}
\label{large-sep}  Let $1\le p\le 2$. Let  $j\ge 4$, $\lambda^{-2/3}\lesssim \mu,{\tilde\mu}\le  c_\circ$, and  $(20d)^{-1}\le c\le (10d)^{-1}$. 
Suppose $\dist (A_k, {\tilde A}_l)\ge c$ and $\dist (A_k, -{\tilde A}_l)\ge c$. Then,  
\Be\label{improved}
\| \chi_{A_k} \fP_\lambda[\psi_j^\kappa ] \chi_{{\tilde A}_l}f\|_{p'}\lesssim   \lambda^{-N} 2^{-Nj}  (\mm)^N \|f\|_p, \quad \forall N>0
\Ee
 for  $\kappa=0,\, \pm, \, \pm \pi$. Moreover,    $(i)$ if $\dist (A_k, {\tilde A}_l)< c$, we  have  \eqref{improved} for $\kappa=\pm \pi,\, 0${\rm ;}  $(ii)$ if  $\dist (A_k, -{\tilde A}_l)< c$,  we  have  \eqref{improved} for $\kappa=\pm,\, 0.$
\end{lem}

To show Lemma \ref{large-sep}, we use the following  elementary lemmas.    

\begin{lem}
\label{trivial}   Let $E$, $F$ be measurable sets.
If  $|T(x,y)|\le D\,\chi_E(x)\chi_F(y)$ for a constant $D$, then $\|T\|_{p\to q}\le D |E|^{1/q} |F|^{1-1/p}$ for $1\le p,  q\le \infty$.
\end{lem}

\begin{lem}\label{S3-osc} Let $N\in \mathbb N_0$, $0<\mu\le 1$, and $L, \lambda>0$. Suppose $a$ is a smooth function supported in an interval $I$ of length $\sim \mu$ and $\phi$ is smooth on $I$. 
If  $|\phi'|\gtrsim   L$, $|(d/ds)^{k+1} \phi |\lesssim  L \mu^{-k}$, and 
$|(d/ds)^k a| \lesssim   \mu^{-k}$,  $0\le k\le  N$  on $I$, then 
\Be
\label{high-order}
\Big|\int e^{i\lambda\phi(s)} a(s)ds\Big|
    \le{C} \mu(1+\lambda\mu L)^{-N}
\Ee 
for a constant $C=C(N)$ independent of  $\lambda, L$ and $\mu$. 
\end{lem}

One can show Lemma \ref{S3-osc}  by repeated integration by parts. It can also be shown by changing  variables $s\to \mu s+s_0$ where $s_0\in \supp a$. Setting  
$\tilde \phi(s)=(\mu L)^{-1}\phi(\mu s+s_0)$ and $\tilde a(s)=a(\mu s+s_0)$, we see    $|\tilde\phi'|\gtrsim  1$, $|(d/ds)^{k} \tilde a| \lesssim 1$, and $|(d/ds)^{k+1} \tilde \phi |\lesssim 1$, $0\le k\le  N,$  on  $\supp\tilde a$. Since the integral equals  $\mu\int e^{i\lambda\mu L\tilde \phi(s)}\tilde a(s)ds$, routine integration by parts yields  \eqref{high-order}.

\begin{proof}[Proof of Lemma \ref{large-sep}]  
We first prove \eqref{improved} for $ \kappa=0,\, \pm, \, \pm \pi$ when $\dist (A_k, {\tilde A}_l)$ $\ge c$ and $\dist (A_k, -{\tilde A}_l)\ge c$.   

Let $(x,y)\in A_k\times {\tilde A}_l$. Note $|x|^2|y|^2\sin^2\!\!\measuredangle(x,y)$ $\ge 2^{-1} c^2$ and $0<\mu,{\tilde\mu}\le 2^{-4}c^2$. By \eqref{angle}  and  \eqref{q-def}, $-\cD(x,y)\sim c^2$ 
and $\cQ(x,y, \tau)\gtrsim c^2$ for $\tau\in \mathbb R$.  By \eqref{ph-d}  we have  $|\partial_s \cP |\gtrsim 2^{2j}$ and also $|\partial_s^l \cP |\lesssim 2^{(l+1)j}$ on $\supp\psi_j^\kappa$. 
Note $ (d/ds)^m (\mathfrak a\psi_j^\kappa)=O(2^{(d/2+m)j})$.  Thus, using Lemma \ref{S3-osc},  we get 
\Be
\label{rapid}
 |\fP_\lambda [\psi_j^\kappa ](x,y)|\lesssim 2^{(d/2-1) j}\big( \lambda2^{j} +1\big)^{-M}, \qquad (x,y)\in (A_k, {\tilde A_l})
\Ee
 for any $M$.  Since $\lambda^{-2/3}\lesssim \mu,{\tilde\mu}\le  c_\circ$,  \eqref{improved} follows by  Lemma \ref{trivial}. 

To prove $(i)$ and $(ii)$,  
    it is sufficient to show $(ii)$ only thanks to \eqref{symmetric2} and the change of variables $y\to -y$. Let $(x,y)\in A_k\times {\tilde A}_l$.   
Since $\dist (A_k, -{\tilde A}_l)$ $< c$,  $|1+\inp xy |\le 3c$  and $|\cD(x,y)|\le 2c$ (by \eqref{angle}). Note $\cos s \ge - \cos 2^{-3}$ on 
$\supp \psi_j^{\kappa }$, $\kappa=\pm, 0$.  Thus, using \eqref{q-def}, we have $|\cQ(x,y, \cos s)|\sim 1$ and $|\partial_s \cP(x,y,s)|\gtrsim  2^{2j}$. 
As before, we also have $|\partial_s^l \cP |\lesssim 2^{(l+1)j}$ for $l\ge 1$  if $s\in \supp\psi_j^{\kappa}$, $\kappa=\pm,0$.  Hence, Lemma  \ref{S3-osc} yields  \eqref{rapid}. 
Consequently,  the estimate \eqref{improved}  follows in the same manner as above.    
\end{proof}

The bounds in Lemma \ref{large-sep} are much smaller  than what  we need  to obtain  for 
$\fP_\lambda$. Recalling  \eqref{coarse-decomp} and \eqref{decomp-proj2},  and discarding the harmless small contributions,  
we need only to consider $ \chi_{A} \fP_\lambda[\psi_j^\kappa ] \chi_{{\tilde A}}$ when $A\subset A_\mu^+$,  ${\tilde A}\subset A_{{\tilde\mu}}^{+}$ are of  diameter $c$ and 
\[
\dist (A, {\tilde A})\le c,  \kappa=\pm, \  \  \text{ or } \dist (A, -{\tilde A})\le c, \kappa=\pm \pi.
\]
By \eqref{symmetric2} and changing  variables $y\to -y$,  the estimate 
for the second case can be deduced from that for the first, so it suffices  to consider the first case only. Moreover, 
the estimate for $ \chi_{A} \fP_\lambda[\psi_j^- ] \chi_{{\tilde A}}$ follows from that for $ \chi_{A} \fP_\lambda[\psi_j^+ ] \chi_{{\tilde A}}$ thanks to  \eqref{symmetric}.
Therefore, the matter is reduced to showing 
 \Be 
 \label{reduced}
\| \sum_{j\ge 4}\chi_{A} \fP_\lambda [\psi_j]  \chi_{{\tilde A}} \|_{p\to p'} \lesssim  
\bpqp \mom^c
\Ee
when $A\subset A_\mu^+$,  ${\tilde A}\subset A_{{\tilde\mu}}^{+}$, and $\dist (A, {\tilde A})\le c$. Henceforth,  we denote  $\psi_j=\psi_j^+$. 
By \eqref{rot} and \eqref{S3-ptkernel},  it follows that 
\Be\label{rot2} \fP_\lambda[\psi_j] (x,y)=\fP_\lambda[\psi_j] (\mathbf U x,\mathbf U y), \quad \mathbf U\in \mathrm O(d). \Ee   
Set $\mathbb S=\{ e'\in  \mathbb S^{d-1}  :  |e' -e_1| < 1/(25d)\}.$ By rotation  we may  also assume 
\[
A, {\tilde A}\subset   \mathbb A_0:=\big \{ x:   |x|^{-1} x \in   \mathbb S \big
\}. 
\]

\subsection{Sectorial decomposition of annuli}
We decompose  $A\times {\tilde A}$ in a way that  we can conveniently 
control the angle between $x$ and $y$.  To do so,  we use a Whitney type decomposition of   $\mathbb S\times \mathbb S$ away from its diagonal. 

Following the typical dyadic decomposition process, for each $\nu\ge 0$ we partition
$\mathbb S$ into spherical caps $\Theta_k^\nu$  such that  $\Theta_k^\nu \subset \Theta_{k'}^{\nu'}$ for some $k'$ whenever $\nu\ge \nu'$ 
and $c_d  2^{-\nu}\le 
  \diam (\Theta_k^\nu)  \le C_d2^{-\nu}$ for some constants $c_d$, $C_d>0$. Let  $\nu_\circ=\nu_\circ(\mu, {\tilde\mu})$ denote 
 the integer  $\nu_\circ$ such that 
\[\mm/2<  2^6 C_d^2 2^{-2\nu_\circ} \le \mm.  
\]
Then, we can write
\[\mathbb S\times \mathbb S=\bigcup_{ \nu:  2^{-\nu_\circ} \le    2^{-\nu} \lesssim 1 }\,\, \bigcup_{k\sim_\nu k'} \Theta_k^\nu\times \Theta_{k'}^{\nu},\] 
where  $k\sim_\nu k'$ means $\dist (\Theta_k^{\nu}, \Theta_{k'}^{\nu})\sim 2^{-\nu}$ if  $\nu> \nu_\circ$ and  $\dist (\Theta_k^{\nu}, \Theta_{k'}^{\nu})\lesssim 2^{-\nu}$ if  $\nu= \nu_\circ $ 
(e.g., see \cite[p.971]{TVV}).  The sets $\Theta_k^{\nu_\circ}$ and $\Theta_{k'}^{\nu_\circ}$, $k\sim_{\nu_\circ} k'$  are not necessarily distanced from each other since  the decomposition 
process terminates  at $\nu=\nu_\circ$.   
Then, it follows that 
\[
A \times  {\tilde A}\,  \subset\!\!\!    \bigcup_{ \nu:   2^{-\nu_\circ}\le     2^{-\nu} \lesssim 1 }\,\, \bigcup_{k\sim_\nu k'} A_k^{\nu}\times \tilde A_{k'}^\nu,
\]
where
\[  A_{k}^{\nu}=\big\{ x\in A_{\mu}^+ :   |x|^{-1} x \in    \Theta_{k}^\nu\big\},   
\quad \  \  \tilde A_{k}^{\nu}=\big\{ x\in A_{{\tilde\mu}}^+:   |x|^{-1} x  \in    \Theta_{k}^\nu\big\}.\] 

Let  $\chi_{k}^{\nu}=\chi_{A_{k}^{\nu}}$ and $\tilde \chi_{k}^{\nu}=\chi_{\tilde A_{k}^{\nu}}.$ The  estimate \eqref{reduced} follows once we obtain 
\Be
\label{eq:la}
\Big\| \sum_{j\ge 4}\sum_{ \nu\le \nu_\circ}\sum_{k\sim_\nu k'} \chi_{k}^{\nu}\, \fP_\lambda [\psi_j] \tilde \chi_{k'}^{\nu} \Big \|_{p\to p'} \lesssim  
\bpqp\mom^c,
\Ee
which we prove in Section  \ref{sec:endpoint}, for $0<\delta(p,p')<\min(2/(d-1), 2/3)$. 

We occasionally use the next elementary  lemma. 

\begin{lem}\label{kkpsum}
    Let  $1\le p\le q\le \infty$. Suppose $\|\chi_{k}^{ \nu}\, \fP_\lambda [\psi_j] \,\tilde \chi_{k'}^{\nu}\|_{p\to q}\le B$  holds whenever $k\sim_\nu k'$. Then for a constant  $C$,  
  \[
\textstyle \|\sum_{k\sim_\nu k'} \chi_{k}^{\nu}\, \fP_\lambda [\psi_j]\, \tilde \chi_{k'}^{\nu} \|_{p\to q}\le C B.
\]
\end{lem}

\begin{proof}  
For each $k$, there are as many as $O(1)$  $k'$ such that  $k\sim_\nu k'$. Thus, 
$ \|  \sum_{k\sim_\nu k'} \chi_{k}^{\nu}\, \fP_\lambda [\psi_j]  \, \tilde\chi_{k'}^{\nu}   f\|_q^q\lesssim \sum_{k\sim_\nu k'} 
 \|   \chi_{k}^{\nu}\, \fP_\lambda [\psi_j] \, \tilde\chi_{k'}^{\nu}   f\|_q^q $. So, it follows that
\[
\textstyle \|  \sum_{k\sim_\nu k'} \chi_{k}^{\nu}\, \fP_\lambda [\psi_j]  \, \tilde\chi_{k'}^{\kappa, \nu} f \|_q^q \lesssim  B^q \sum_{k'}  \| \tilde\chi_{k'}^{\nu} f\|_p^q.
\]
The right hand side is clearly bounded by $CB^q\|f\|_p^q$ since $p\le q$ and $\{\tilde A_k^{\nu}\}_k$  are disjoint. 
\end{proof}

\subsection{The kernel of $\mathfrak P_\lambda[\psi_j]$}
In this subsection, we obtain estimates for the kernel $\mathfrak P_\lambda[\psi_j](x,y)$, which we use later.

\begin{lem}
\label{case-a}  Let   $0<{\tilde\mu}\le \mu\ll 2^{-\nu}\le 1/100$, and  $(x,y)\in  A_k^{\nu}\times {\tilde A}_{k'}^{\nu}$, $k\sim_\nu k'$. Then for any $N>0$, we have  
\begin{equation}
\label{a-case-kernel}
   |\mathfrak P_\lambda[\psi_j](x,y)|    \lesssim  2^{\frac{d-2}{2}j} \big( \lambda2^{j}\max(2^{-2\nu}, 2^{-4j}) +1\big)^{-N} .
                         \end{equation}
\end{lem}

\begin{proof} Note $\measuredangle(x,y)\sim 2^{-\nu}$  for $(x,y)\in  A_k^{\nu}\times {\tilde A}_{k'}^{\nu}$, $k\sim_\nu k'$.   Since $\mu\ll 2^{-\nu}$, it is easy to see that $|x-y|\sim 2^{-\nu}$. So, 
 $-\cD(x,y)\sim 2^{-2\nu}$ by \eqref{angle}.  Note  $|x-y|^2+\cD(x,y)=(1-\inp xy)^2$, hence  $|1-\inp{x}{y}|\lesssim 2^{-\nu}.$
Combining these observations with  \eqref{ph-d} and \eqref{q-def},  we have
\[
|\partial_s \cP(x,y,s) |
            \gtrsim  2^{2j} \max(2^{-2\nu}, 2^{-4j}), \quad s\in \supp \psi_j. 
\]
By \eqref{ph-d},  it also follows that   
$| \ddk(x,y,s)| \lesssim  2^{(1+k)j}  \max(2^{-2\nu}, 2^{-4j})$ for $s\in \supp \psi_j$. 
Thus,  using  Lemma  \ref{S3-osc} with $L=2^{2j} \max(2^{-2\nu}, 2^{-4j})$ and $\mu=2^{-j}$, we get  \eqref{a-case-kernel} in the same manner as in the proof of Lemma \ref{large-sep}. 
\end{proof}

\begin{lem} 
 Let 
$0< \tilde\mu\ll \mu\le 1/100$ and $(x,y)\in A_k^{ \nu}\times \tilde A_{k'}^{\nu}$, $k\sim_\nu k'$. Suppose $2^{-\nu_\circ}\le 2^{-\nu}\lesssim \mu$. Then for any $N>0$,  we have
    \begin{equation}
\label{b-case-kernel}
     |\mathfrak P_\lambda[\psi_j](x,y)|      \lesssim  \begin{cases}   
                                 2^{\frac{d-2}{2}j} \big( \lambda 2^{j}\mu^{2}+1\big)^{-N} ,        &  \  \  2^{-j}  \ll   \smu  \ \ 
                                                                \\[2pt]
                                   2^{\frac{d-2}{2}j} \big( \lambda 2^{-3j}+1\big)^{-N},        & \  \  2^{-j}\,  \gg    \smu
                          \end{cases}. 
\end{equation} 
 Additionally, if $2^{-j}\sim \sqrt{\mu}$, and  $\mathcal D(x,y)\sim\mu\tilde\mu$ or $\mathcal D(x,y) <0$, then
    \Be\label{b-case-sqrtmu}
    |\mathfrak P_\lambda[\psi_j](x,y)|\lesssim \lambda^{-\frac12}\mu^{\frac{1-d}{4}}|\mathcal D(x,y)|^{-\frac14}. 
    \Ee
\end{lem}

\begin{proof} 
We  consider \eqref{b-case-kernel} first.  To this end, 
we claim  
\Be
\label{lower-case}
|\partial_s \cP(x,y,s) |
            \gtrsim   \begin{cases}   
                                 \mu^2 2^{2j},            &  \   2^{-j}   \ll  \smu  
                                                                 \\
                                  \   2^{-2j},         &  \  2^{-j}\,   \gg   \smu
                          \end{cases},  \quad \  s\in \supp\psi_j. 
  \Ee
 Note  $2(1-\inp xy) = 1-|x|^2+1-|y|^2 +|x-y|^2$. Since $|x-y|\sim \mu,$ $|1-\inp xy|\sim \mu$. So,  
$
|\partial_\tau \mathcal Q(x,y,\tau)| = 2|\tau-\inp xy|\lesssim \mu
$
if $\tau\in [1-c\mu,1]$ for some $c>0$. If $ 2^{-j}\ll \sqrt{\mu}$, by the mean value theorem 
$| \cQ(x,y,\cos s)- \cQ(x,y,1) |\ll \mu^2 $ for  $s\in\supp{\psi_j}$  because $|1-\cos s|\le s^2/2$.  Observing   
$\cQ(x,y,1) = |x-y|^2 \sim \mu^2$, we see   $\cQ(x,y,\cos s)\sim \mu^2$  for  $s\in\supp{\psi_j}$ if $ 2^{-j}\ll \sqrt{\mu}$. Thus, by \eqref{ph-d} we have 
the first case  in \eqref{lower-case}.  For the second case, note $|\mathcal D(x,y)|\lesssim 2^{-2\nu}\lesssim \mu^2$  and $1-\tau \sim 2^{-2j}$  for $\tau\in\cos{(\supp{\psi_j})}$. Recalling \eqref{q-def}, we have 
 $\cQ(x,y,\tau)=(\tau-\inp xy)^2-\cD(x,y) \sim 2^{-4j}$ if $2^{-j}\gg \sqrt{\mu}$.  So, by \eqref{ph-d} we obtain the second case  in \eqref{lower-case}.

A computation using \eqref{ph-d} shows
\Be 
\label{lower-case2}
\Big| \ddk(x,y,s)\Big| \lesssim   \begin{cases}   
                                 \mu^2 2^{(1+k)j} ,           
                                          &  \  2^{-j}   \ll  \smu  
                                         \\[2pt]
                                \ 2^{-4j}  2^{(1+k)j},         
                                &  \  2^{-j}   \gg  \smu  
                          \end{cases} ,  \quad  \  \  s\in \supp\psi_j.   \Ee
Therefore, combining \eqref{lower-case} and \eqref{lower-case2},  we obtain \eqref{b-case-kernel} by Lemma \ref{S3-osc}. 

We now turn to \eqref{b-case-sqrtmu} and consider the case $\mathcal D(x,y)\sim \mu\tilde\mu$ first. Since $|1-\inp xy|\sim \mu$,  $\cQ(x,y,\cdot)$ has two distinct zeros  $r_1> r_2$, which are close to $1$. 
Let $s_1$, $s_2\in (0, \pi/2)$ be numbers such that $\cos s_i =r_i$, $i=1,2$,  and $s_1< s_2$.  By \eqref{ph-d} and \eqref{q-def} we have 
\Be
\label{eq:phase} 
 \partial_s \mathcal P(x,y,s)= -\frac{\int_{s_1}^s \sin \tau d\tau \int_{s_2}^s \sin \tau d\tau}{2\sin^2s}.
 \Ee

We decompose the integral $I_j:=\mathfrak P_\lambda[\psi_j](x,y)$ away from  $s_1, s_2$.  Let us set  
$  \psi^{k, \pm}_l(s)= \psi(\pm 2^k(s-s_l)),$ $l=1,2$. Let $c>0$ be a small constant
and $k_0$ be the smallest integer such that $2^{-k_0}\le c\tmu^{1/2}$. 
Then, we  put  \[ \psi_\ast=1-\ssum{k} \psi^{k,-}_1  -  \ssum{k\ge k_0} \psi^{k,+}_1-\ssum{k\ge k_0} \psi^{k,-}_2 - \ssum{k} \psi^{k,+}_2,\] so 
$\supp \psi_\ast\subset [s_1+c_1 \tmu^{1/2}, s_2-c_1 \tmu^{1/2}]$ for a constant $c_1>0$. Also  set  
 \begin{align*} 
{I}_{l}^{k, \pm}=\int (\psi^{k, \pm}_l \psi_j\mathfrak a)(s)e^{i\lambda\mathcal P(x,y,s)}ds, \quad  \  \ {I}_{\ast} =\int (\psi_\ast\psi_j\mathfrak a)(s)e^{i\lambda\mathcal P(x,y,s)}ds. 
\end{align*}
Then,  $ I_j=I_{\ast}+I_{1}^-+I_{1}^+ + I_{2}^-+I_{2}^+$, where  
 \[   I_{1}^-=  \ssum{k} {I}_{1}^{k,-},   \quad I_{1}^+= \ssum{k\ge k_0} {I}_{1}^{k,+},  \quad  I_{2}^-=\ssum{k\ge k_0} {I}_{2}^{k,-},  \quad I_{2}^+= \ssum{k} {I}_{2}^{k,+}.\]

 Since $\cos s_1-\cos s_2=2\sqrt \cD\sim (\mu\tmu)^{1/2}$ and $s_1,s_2\sim \sqrt \mu\sim 2^{-j},$  $s_2-s_1\sim \tmu^{1/2}$.  By \eqref{eq:phase} it follows that $ |\partial_s \mathcal P(x,y,s)|\gtrsim   \tmu $ for  $s\in \supp \psi_\ast$. Since $\|  (\psi_\ast  \mathfrak a  \psi_j)'  \|_1$ $\lesssim 2^{dj/2}$,   
van der Corput's lemma (e.g., \cite[p. 334]{St93}) gives $|I_{\ast}|\lesssim  \min( \lambda^{-1}2^{(d-4)j/2} \mathcal D^{-1},$ $2^{(d+2)j/2} \cD^{1/2})$, which yields   $|I_{\ast}|\lesssim \lambda^{-1/2}\mu^{(1-d)/4} \mathcal D^{-1/4}$. Therefore, we have only  to show
 \[ |I_{l}^\pm|\lesssim \lambda^{-1/2}2^{(d-1)j/2} \mathcal D^{-1/4},  \  \quad  l=1,2.\]
 We  consider ${I}_{1}^{-}$ only. The estimates for the others can be shown  in a similar manner.  Since $s_2-s_1\sim \tmu^{1/2}$, by \eqref{eq:phase}  we have
$
 |\partial_s \mathcal P(x,y,s)|
\gtrsim  2^j2^{-k}\cD^{1/2}
 $
 for $s\in \supp \psi^{k,-}_1$.    
By van der Corput's  lemma, $  |{I}_{1}^{k,-} |\lesssim    2^{dj/2}\min( \lambda^{-1} 2^{k-j} \mathcal D^{-1/2}, $ $ 2^{-k} ). $ 
Summation over $k$ gives $|{I}_{1}^{-}|\lesssim  \lambda^{-1/2}2^{(d-1)j/2}$ $\mathcal D^{-1/4}$ as desired.

Following  the previous argument closely, we  show \eqref{b-case-sqrtmu} when $\cD<0$. 
From \eqref{ph-d} and \eqref{q-def}, we have
\Be\label{ggps}
|\partial_s\mathcal P(x,y,s)|\sim ((\cos s - \inp xy)^2+|\cD|)\mu^{-1}, \quad s\in  \supp{\psi_j}. 
\Ee
Let  $s_*\in (0,\pi/2)$ denote the point such that  $\cos s_* = \inp xy$, and let $k_*$  be the smallest number satisfying  $2^{-k_*}\le  |\cD|^{1/2}\mu^{-1/2}$.
We decompose the integral $I_j$, using the cutoff functions $
\psi^k(s) = \psi(2^k(s-s_*)) + \psi(2^k(s_*-s)),$ $k< k_\ast,$ 
and $\psi^{k_\ast}(s) := 1-\ssum{k\le k_*}\psi^k(s)$.  As before, setting
\[ I^k := \int (\psi^k \psi_j\mathfrak a)(s)e^{i\lambda\mathcal P(x,y,s)}ds,\quad I^{k_\ast}:=\int (\psi^{k_\ast} \psi_j \mathfrak a)(s)e^{i\lambda\mathcal P(x,y,s)}ds, \]
we break $I_j=\sum_{ k\le k_\ast } I^k$. From \eqref{ggps}  we see $|\partial_s\mathcal P(x,y,s)|\gtrsim 2^{-2k}$ for $s\in\supp{\psi^k}$. Thus, the van der Corput lemma gives $|I^k|\lesssim \mu^{{-d}/{4}}\min(\lambda^{-1}2^{2k}, 2^{-k})$ $\le \mu^{{-d}/{4}}\lambda^{-1/2}2^{k/2}$. 
Taking  sum over $k\le k_\ast$, we get \eqref{b-case-sqrtmu}. 
\end{proof}

\subsection{An $L^2$-estimate for $\fP_\lambda[ \eta_j ]$} We denote  $A_{\mu}^\circ=\{x:   |1- |x|| \le  2\mu\}$ and  $A_{\lambda, \mu}^\circ=\{x: \lambda^{-1/2}x\in  A_{\mu}^\circ\}$. We also set  
\[ \chi_\mu^\circ=\chi_{A_\mu^\circ},  \qquad \chi_{\lambda, \mu}^\circ=\chi_{A_{\lambda, \mu}^\circ}.\]

 \begin{lem}
 \label{l2-mmj23}   Let $\lambda^{-\frac23}\lesssim\tilde\mu,  \mu\le 1/4$ and $2^{-j}\gtrsim (\lambda\mu)^{-1}$. 
 Suppose $\eta_j \in C_c^\infty (-\pi, \pi)$ supported in an interval of length $\sim 2^{-j}$ satisfies $|\eta_j^{(k)}|\lesssim 2^{jk}$ for any $k$. 
 Then, 
\begin{equation}\label{jl222-}    \| \chi_{\mu}^\circ \mathfrak P_\lambda [ \eta_j ]  \chi_{\tilde\mu}^\circ   \|_{2\to 2} \lesssim   \lambda^{-\frac d2}\mmf    .   
 \end{equation}
\end{lem}

To prove \eqref{jl222-}, we instead show 
$ \|   \chi_{\lambda, \mu}^\circ \Pi_\lambda [ \eta_j ]  \chi_{\lambda, \tilde\mu}^\circ  \|_{2\to 2} \lesssim   (\mu{\tilde\mu})^{1/4} $      
which  is equivalent to  \eqref{jl222-} (see \eqref{eq:norm-scaled}). In fact, we can show  a stronger estimate 
\Be
\label{stronger}   
\|\chi_{\lambda, \mu}^e \Pi_\lambda [ \eta_j ]  \chi_{\lambda, \tilde\mu}^e\|_{2\to2} \lesssim  (\mu\tilde\mu)^\frac14
\Ee
for $\lambda^{-\frac23}\lesssim\tilde\mu,  \mu\le 1/4$, where
 \[  
 \chi_{\lambda,\mu}^e=\chi_{A_{\lambda,\mu}^e}, \quad \ \ A_{\lambda,\mu}^e =\{x:|x|\ge\lambda^{1/2}(1-\mu)\}.
 \]

We now recall the next estimates which follow from \cite[Theorem 3]{T05}:  
\begin{align}
\label{smallmu} 
& \, \, \|\chi_{\lambda,  \lambda^{-2/3} }^\circ \Pi_\lambda \chi_{\lambda, \lambda^{-2/3}}^\circ\|_{p\rightarrow p'}\lesssim  \lambda^{\frac{d}{6}\delta(p,p')-\frac 13},
\\
\label{est:ext--}
    \|\chi_{\lambda, \mu}^-\Pi_\lambda &\chi_{\lambda, \mu}^-\|_{p\rightarrow p'}\lesssim \lambda^{\frac d6\delta(p,p')-\frac13}\big(\lambda^{2/3}\mu (1+\mu) \big)^{-M},  \  \lambda^{-2/3}\lesssim \mu
\end{align}
for $1\le p\le 2$. The estimate \eqref{est:ext--} holds for any  $M>0$.  In particular,  the $L^2$--$L^\infty$ estimate in  \cite[Theorem 3]{T05} implies \eqref{est:ext--} for $p=1$. Interpolation with 
$L^2$ estimate shows \eqref{est:ext--} for $1\le p\le 2$.  One can also show \eqref{smallmu}  and \eqref{est:ext--} in an elementary manner   using  estimates for the kernels of $\Pi_\lambda$   (e.g., see  \cite{JLR-hermite}).

\begin{proof}[Proof of  \eqref{stronger}] Considering the adjoint operator, we may assume $\tilde\mu\le \mu$.   We begin by showing
\Be \label{stronger1}  \|  \chi_{\lambda, \mu}^e \Pi_\lambda  \chi_{\lambda, {\tilde\mu}}^e \|_{2\to 2} \lesssim (\mu{\tilde\mu})^\frac14, \qquad     \lambda^{-2/3}\lesssim \tilde\mu\le \mu\le 1.\Ee
Note $ \|  \chi_{\lambda, \mu}^e \Pi_\lambda  \chi_{\lambda, {\tilde\mu}}^e \|_{2\to 2}\le 
\|  \chi_{\lambda, \mu}^e \Pi_\lambda \|_{2\to 2} \|   \Pi_\lambda \chi_{\lambda, {\tilde\mu}}^e \|_{2\to 2}$. By duality it is enough to show  $\|\chi_{\lambda, \mu}^e \Pi_\lambda \|_{2\to 2}\lesssim \mu^{1/4}$ for $\lambda^{-2/3}\lesssim \mu\le 1$. 
To do this, recalling that $\mu, \tilde\mu$ denote dyadic numbers,   we decompose 
 \[ \chi_{\lambda, \mu}^e = \ssum{\lambda^{-2/3}\le{\tilde\mu} \le \mu} 
                \chi^+_{\lambda, {\tilde\mu}}  + \ \chi_{\lambda,\lambda^{-2/3}}^\circ+  \ssum{\lambda^{-2/3}\le {\tilde\mu}< 1 } \chi^-_{\lambda, {\tilde\mu}} \ + \  \  \ssum{1\le {\tilde\mu}  }  \chi^-_{\lambda, {\tilde\mu}}.\]
By \eqref{est:ext--}, it follows that $ \sum_{1\le {\tilde\mu}  } \|     \chi^-_{\lambda, {\tilde\mu}} \Pi_\lambda \|_{2\to 2}\lesssim \lambda^{-N}$ for any $N$. 
The estimates \eqref{kt-nu},  \eqref{smallmu}, and  \eqref{est:ext--} yield   $\sum_{\lambda^{-2/3}\le{\tilde\mu} \le \mu} 
             \|   \chi^+_{\lambda, {\tilde\mu}}   \Pi_\lambda\|_{2\to 2}\lesssim  \mu^{1/4}$, $  \|   \chi_{\lambda,\lambda^{-2/3}}^\circ   \Pi_\lambda\|_{2\to 2}\lesssim  \lambda^{-1/6}$, 
             and $\sum_{\lambda^{-2/3}\le {\tilde\mu} \le 1}$ $\| \chi^-_{\lambda, {\tilde\mu}} \Pi_\lambda \|_{2\to2}\lesssim \lambda^{-1/6}$, respectively.  
             Therefore,  we get the desired estimate  
                         since $ \mu\gtrsim \lambda^{-2/3}$.

Now, observe $
\Pi_\lambda [\eta_j]  =  \sum_{\lambda'} \widehat {\eta_j} ( 2^{-1}(\lambda'-\lambda))\Pi_{\lambda'},
$
which follows by  \eqref{proj-pi} and \eqref{schro-op}. Since $|\widehat {\eta_j}(\tau)|\lesssim  2^{-j} (1+2^{-j} |\tau|)^{-N}$,  we have
\[ \| \chi_{\lambda, \mu}^e \Pi_\lambda[\eta_j]  \chi_{\lambda, {\tilde\mu}}^e\|_{2\to 2}\lesssim   2^{-j} \sum{}_{{}_{\mathlarger \lambda' }} \left(1+2^{-j}|\lambda-\lambda'|\right)^{-N}\|\mmfp\|_{2\to 2} .
\]  
However, we can not directly apply \eqref{stronger1}  since $\lambda\neq \lambda'$.
We get around the problem  by enlarging the sets $A_{\lambda,\mu}^e, A_{\lambda,\tilde\mu}^e$.
Let
\[
\ell_{\lambda'}(\rho)=
\begin{cases}
   \  \  \   (\lambda/\lambda')^{\frac23}\rho, \quad &\lambda \ge \lambda',
    \\[2pt]
   \,   (\lambda'-\lambda)/\lambda'+\rho, \quad &\lambda < \lambda'.
    \end{cases} \]
     Note  that  $\ell_{\lambda'}(\rho) \gtrsim(\lambda')^{-2/3}$ for $\lambda^{-2/3}\lesssim\rho$. 
 Since  $(\lambda')^{1/2}(1-\ell_{\lambda'}(\mu))\le \lambda^{1/2}(1-\mu)$, 
i.e., $ A_{\lambda,\mu}^e\subset A_{\lambda',\ell_{\lambda'}(\mu)}^e,$  $\|\mmfp\|_{2\to 2}
\le\|\chi^e_{\lambda',\ell_{\lambda'}(\mu)}\Pi_{\lambda'}\chi^e_{\lambda',\ell_{\lambda'}(\tilde \mu)}\|_{2\to 2}$.  
Using \eqref{stronger1}, we have $ \|\mmfp\|_{2\to 2}   \lesssim    ( \ell_{\lambda'}(\mu)\ell_{\lambda'}(\tilde \mu))^{1/4}$.  
 Therefore,   it suffices for  \eqref{stronger} to show 
\[   2^{-j} \ssum{\lambda'} \left(1+ 2^{-j} |\lambda-\lambda'|\right)^{-N}  \big( \ell_{\lambda'}(\mu)\ell_{\lambda'}(\tilde \mu)\big)^{\frac14}\lesssim  (\mu\tilde\mu)^{\frac14}.\]       
This can be shown  by a simple computation because $2^{-j}\gtrsim (\lambda\mu)^{-1}$ and $\mu,{\tilde\mu}\gtrsim\lambda^{-2/3}$. We omit the detail. 
\end{proof}

\subsection{Proof of  Corollary \ref{cor:weighted}}  
\label{pfpf}
Before we conclude this section,   assuming   Theorem \ref{unb},  we prove   Corollary \ref{cor:weighted}.  
We follow the lines of argument for the endpoint bound which is sketched in the introduction.

\begin{proof}[Proof of Corollary \ref{cor:weighted}]
We  first  prove \eqref{weighted} for $2< q\le q_\ast:=2(d+1)/(d-1)$. 
Let  us set
\[\mathcal W=w_+^{\frac{d+3}4\delta(2,q)-\frac14}.\]
Note that $\mathcal W= w_+^\gamma$ 
if $2< q\le q_\ast$.  Recall $  \mu_\circ =\lambda^{-2/3}$ and $c_\circ=  (100d)^{-2}.$ By the triangle inequality, $\|\mathcal W w_{-}^N \Pi_\lambda\|_{2\to q}$ is bounded by 
\[ 
\|  \sum_{\mu\ge \mu_\circ} \mathcal W \chi_{\lambda,\mu}^+ \Pi_\lambda\|_{2\to q}+ \| \chi_{\lambda, \mu_\circ}^\circ \Pi_\lambda \|_{2\to q} 
+ \sum_{2^k\ge \mu_\circ}  ( \lambda^\frac232^k(2^k+1))^N \| \chi_{\lambda, 2^k}^- \Pi_\lambda \|_{2\to q},\]
where $\mu\in \mathbb D^{-}$. By the estimates \eqref{smallmu} and \eqref{est:ext--} the last two are bounded by $C \lambda^{(d\delta(2,q)-1)/6}$.    
It follows by  \eqref{kt-nu} that $  \sum_{\mu\ge  c_\circ}\|  \mathcal W \chi_{\lambda,\mu}^+ \Pi_\lambda\|_{2\to q}\lesssim \lambda^{(d\delta(2,q)-1)/6}$.
Therefore,  the matter is reduced to showing 
\[
  \| \mathfrak W\|_{2\to q}\lesssim \lambda^{\frac d6\delta(2,q)- \frac 16}, 
\]
where 
$ \mathfrak W=\sum_{ \mu_\circ< \mu < c_\circ} \mathcal W \chi_{\lambda,\mu}^+ \Pi_\lambda.$
Equivalently, we need to   show 
\Be 
\label{weighted0}  
\|\mathfrak W\mathfrak W^\ast \|_{q' \to q} \lesssim  \lambda^{\frac d6\delta(q',q)-\frac13}.
\Ee

To this end, we write 
 \[ 
\textstyle \mathfrak W\mathfrak W^\ast =  \sum_{k} \mathfrak  (\mathfrak W\mathfrak W^\ast)_k:=
 \sum_{k} \sum_{(\mu, {\tilde\mu})\in \mathfrak D_k}  \mathcal W \chi_{\lambda,\mu}^+ \Pi_\lambda  \mathcal W \chi_{\lambda,{\tilde\mu}}^+\,.\]
 Recall that
 $\mathfrak D_k=\{(\mu, \tilde\mu):  {{\tilde\mu}}/\mu\in [ 2^k, 2^{k+1}), \mu, \tilde\mu\in (\mu_\circ, c_\circ)\}$ and $\mu, \tilde\mu\in \mathbb D$.
   Since  $\supp \chi_{\lambda,{\tilde\mu}}^+$ are almost disjoint, for each $k$ we have
\[
\textstyle \big\|(\mathfrak W\mathfrak W^\ast)_k   f\big\|_{q}
        \lesssim  \big(\sum_{(\mu, {\tilde\mu})\in \mathfrak D_k} \|  \mathcal W \chi_{\lambda,\mu}^+ \Pi_\lambda  \mathcal W \chi_{\lambda,{\tilde\mu}}^+   f\|_{q}^{q}\big)^{1/{q}} .
        \]
 Note $ \mathcal W \chi_{\lambda,\mu}^+  \sim  \lambda^{\frac {d+3}{12}\delta(q',q)-\frac16} \mu^{\frac{d+3}{8} \delta(q',q)-\frac14} \chi_{\lambda,\mu}^+$. 
Using \eqref{unbalanced}, we get
 \[ \|  \mathcal W \chi_{\lambda,\mu}^+ \Pi_\lambda  \mathcal W \chi_{\lambda,{\tilde\mu}}^+   f\|_{q}  \lesssim \lambda^{\frac d6\delta(q',q)-\frac13} 2^{-c|k|} \|\chi_{\lambda,{\tilde\mu}}^+   f\|_{q'}, 
 \quad {{\tilde\mu}}/\mu \sim 2^k. \] 
Indeed, when $2^k>1$, we  consider  the adjoint operator $ T^\ast=\mathcal W \chi_{\lambda,{\tilde\mu}}^+ \Pi_\lambda  \mathcal W \chi_{\lambda,\mu}^+$ of  $T:=  \mathcal W \chi_{\lambda,\mu}^+ \Pi_\lambda  \mathcal W \chi_{\lambda,{\tilde\mu}}^+ $ and then we may use the estimate \eqref{unbalanced} thanks to the fact that $ \| T\|_{q'\to q}=\|  T^* \|_{q'\to q}$.  Since   $\sum_{ {\tilde\mu}}  \|\chi_{\lambda,{\tilde\mu}}^+   f\|_{q'}^{q} \lesssim \|f\|_{q'}^q$, we obtain 
\begin{align*}
\big\|(\mathfrak W\mathfrak W^\ast)_k   f\big\|_{q} 
       \lesssim   2^{-c|k|} \lambda^{\frac d6\delta(q',q)-\frac13}  \|f\|_{q'}. 
                \end{align*}
 Summation  over $k$ gives  
the desired estimate  \eqref{weighted0}.

We now prove  \eqref{weighted} for $ 2(d+1)/(d-1) <q\le \infty$. 
Note $\gamma= \tfrac12- \tfrac{d}2\delta(2,q).$  By decomposing the operator in the same way as above, it suffices to show  
\Be 
\label{weighted3} 
\textstyle  \|  \sum_{ \mu_\circ < \mu<   c_\circ  } w_+^{\frac12- \frac{d}2\delta(2,q) } \chi_{\lambda,\mu}^+ \Pi_\lambda f\|_{2\to q}\lesssim \lambda^{\frac d6\delta(2,q)-\frac16} \|f\|_{q'}
\Ee
for $2(d+1)/(d-1) <q\le \infty$. The other parts can be handled similarly as before.  By interpolation, we need only to show \eqref{weighted3} for $q=\infty,$ $2(d+1)/(d-1)$.
Thanks to disjointness of the annuli, \eqref{kt-nu} for $q=\infty$ gives \eqref{weighted3} for $q=\infty$, while  \eqref{weighted0}   is equivalent to  \eqref{weighted3} when $q=2(d+1)/(d-1)$. 
\end{proof}

\section{Asymmetric improvement: Proof of Proposition \ref{improvedbound} } 
\label{sec:endpoint}

We prove Proposition  \ref{improvedbound} by establishing the estimate \eqref{eq:la}. To this end,  
we separately consider some cases.   The desired estimates can be shown by  the kernel estimates (in the previous section) except for the case   $2^{-j}\sim \sqrt \mu$,  $2^{-\nu}\sim \mmt$,  and $|\cD|<\vepc \mm$, 
which requires  a different approach. We handle this case  in the next section.

To show \eqref{eq:la},  we  distinguish the cases $\nu\in \nb$ and  $\nu\in \nc$, where 
\[  \nb = \{ \nu: 2^{-\nu}\gg \mmt \text{ or }   \nu=\nu_\circ\}  , \quad  \quad  \nc=   \{ \nu: \mmt \gtrsim  2^{-\nu} > 2^{-\nu_\circ}\}  . \]

\subsection{The sum over  $\nu\in\nb$}  In this case,  the desired estimates are easier. 

\begin{prop} 
\label{lala}  
Let $d\ge 2$ and ${\tilde\mu}, \mu$ satisfy \eqref{mucon}.   If  ${2}/{(d+3)}<\dpp<  {2}/{(d-1)}$, then for some $c>0$,  we have 
\Be 
\label{eq:lala}
\|  \sum_{ \nu\in\nb}\sum_{j\ge 4}\sum_{k\sim_\nu k'} \chi_{k}^{\nu} \fP_\lambda [\psi_j] \tilde\chi_{k'}^{\nu} \|_{p\to p'} \lesssim  
\bpqp
\mom^{c}. 
 \Ee
\end{prop}

\begin{proof}  We  further divide  $\nb=\nb^1\cup \nb^2$, where 
\[ \nb^1=\{ \nu: 2^{-\nu}  \gg   \mu\} , \quad  \  \ \nb^2=\{ \nu:  \mu\gtrsim   2^{-\nu}\gg  \mmt   \text{ or } \nu=\nu_\circ\}.   \] 

We first consider    $\nu\in \nb^1$. 
Using Lemma \ref{case-a},  we get 
\[ 
\|\chi_{k}^{\nu} \mathfrak P_\lambda [ \psi_j ]  \tilde\chi_{k'}^{\nu}\|_{p\to p'} 
          \lesssim  \begin{cases}   
                                 \lambda^{\frac{d-1}2\dpp-\frac d2} 2^{(\frac{d-1}2\dpp -1)j}  2^{\nu\dpp },  
                                      & \  \ \  \!\! 2^{-2j} 
\lesssim   2^{-\nu}                                  
\\[2pt]
                                 \lambda^{\frac{d-1}2\dpp-\frac d2} 2^{(\frac {d+3}2\dpp -1) j},    
                                 &  \  \  \  \!\!  2^{-2j}  \gg 2^{-\nu}
                          \end{cases}, \hspace{-.4cm}
\]
for  $k\sim_\nu k'$ and $1\le p\le 2$.  Indeed,  taking $N=1/2$ in  \eqref{a-case-kernel}, we get the $L^1$--$L^\infty$ bounds, and then interpolation with $\|\chi_{k}^{\nu} \mathfrak P_\lambda [ \psi_j ]  \tilde\chi_{k'}^{\nu}\|_{2\to 2}\lesssim \lambda^{- d/2} 2^{-j}$, which follows from 
\eqref{easy-l2} and \eqref{eq:norm-scaled}, gives the estimates for $1\le p\le 2$.  By Lemma \ref{kkpsum}, we have the same bounds on  $\|\sum_{k\sim_\nu k'}\chi_{k}^{\nu} \mathfrak P_\lambda [ \psi_j]  \tilde\chi_{k'}^{\nu}\|_{p\to p'}$ as above. Since  ${2}/{(d+3)}<\dpp<  {2}/{(d-1)},$  summation over $j$ gives  
\[
\sum_{j}  \| \sum_{k\sim_\nu k'}\chi_{k}^{\nu} \mathfrak P_\lambda [ \psi_j]  \tilde\chi_{k'}^{\nu}\|_{p\to p'}
\le C  \lambda^{\frac{d-1}2\dpp-\frac d2}    2^{(\frac{d+3}4\dpp-\frac12)\nu}. 
\]  
Thus, taking sum over  $\nu: 2^{-\nu}\gg \mu$, we obtain
\Be
\label{nb1}
 \sum_{\nu\in \nb^1}  \sum_j \big\|   \sum_{k\sim_\nu k'}\chi_{k}^{\nu} \mathfrak P_\lambda [ \psi_j]  \tilde\chi_{k'}^{\nu} \big\|_{p\to p'} 
  \lesssim    \bpqp \mom^{\frac{d+3}8\dpp -\frac14}. 
 \Ee 

We turn to the case   $\nu\in \nb^2$. As above, by \eqref{b-case-kernel} and   \eqref{b-case-sqrtmu} we have
\begin{equation}
\label{eq:sec8} 
       \,\, \|\chi_{k}^{\nu} \mathfrak P_\lambda [ \psi_j ]  \tilde\chi_{k'}^{\nu}\|_{p\to p'}          
              \!  \lesssim \! \begin{cases}                                 
                      \! \lambda_p\,  2^{(\frac{d-1}2\dpp -1)j}  \mu^{-\dpp} ,    &  2^{-j}   \ll  \smu            
                                                                \\                                
                       \! \lambda_p\,  \mu^{\frac14-\frac d4\dpp}\tilde \mu^{\frac14-\frac14\dpp}  2^{ \frac \nu 2 \dpp },  &  2^{-j}  \sim  \smu 
                                                                 \\                        
                       \!  \lambda_p\,  2^{(\frac {d+3}2\dpp  -1) j} ,         & 2^{-j}\,   \gg   \smu    
\end{cases} \!\!\! \end{equation}
for $k\sim_\nu k'$ and $1\le p\le 2$ when $\nu\in \nb^2$. Here,  $ \lambda_p$ denotes $\lambda^{\frac{d-1}2\dpp-\frac d2}.$ 
For the first and third cases, we take $N=1/2$ in  \eqref{b-case-kernel} to get  the $L^1$--$L^\infty$ estimates and  interpolate them  with $\|\chi_{k}^{\nu} \mathfrak P_\lambda [ \psi_j ] \tilde  \chi_{k'}^{\nu}\|_{2\to 2}\lesssim \lambda^{-\frac d2} 2^{-j}$. 
We get the second  case using $\|\chi_{k}^{\nu} \mathfrak P_\lambda [ \psi_j ]  \tilde\chi_{k'}^{\nu}\|_{2\to 2} \lesssim  \lambda^{-d/2} (\mu\tmu)^{1/4}$ and 
$\|\chi_{k}^{\nu} \mathfrak P_\lambda [ \psi_j ]  \tilde\chi_{k'}^{\nu}\|_{1\to \infty}$ $\lesssim  \lambda^{-1/2}   \mu^{(1-d)/{4}}2^{\nu/2}$, which follow
from \eqref{jl222-} and  \eqref{b-case-sqrtmu}, respectively.

  Since ${2}/{(d+3)}<\dpp<  {2}/{(d-1)}$,  
 by Lemma \ref{kkpsum}  we get
\[
    \sum_{j}\|\sum_{k\sim_\nu k'} \chi_{k}^{\nu} \mathfrak P_\lambda [ \psi_j] \tilde \chi_{k'}^{\nu} \|_{p\to p'}\lesssim \bpqp  \mathrm B_\ast(\mu, \tilde\mu), \qquad \nu\in \nb^2, 
    \]
where 
\[ \mathrm B_\ast(\mu, \tilde\mu)=(\tilde \mu/\mu)^{\frac{d+3}8\dpp -\frac14} + \mu^{-\frac{d-3}8\dpp}\tilde \mu^{\frac{d+1}8\dpp}2^{\frac{\nu}{2}\dpp}.\] Thus, $\sum_{\nu\in\nb^2}\sum_{j}\|\sum_{k\sim_\nu k'} \chi_{k}^{\nu} \mathfrak P_\lambda [ \psi_j]\tilde\chi_{k'}^{\nu}\|_{p\to p'}$ is bounded by 
\begin{align*}
 C\,    \bpqp 
   \big(\mom^{\frac{d+3}8\dpp -\frac14}\log{(\mu/\tilde\mu)} + \mom^{\frac{d-1}{8}\dpp}\big).
\end{align*}
Combining this and \eqref{nb1}, we obtain \eqref{eq:lala} for some $c>0$ .
\end{proof}

\subsection{The sum over $\nu\in \nc$}  Since there are only as many as $O(1)$ $\nu$,  it suffices  to consider a single $\nu$. 

\begin{prop}
\label{l0} 
Let $d\ge 3$,  $\mu, \wmu$ satisfy \eqref{mucon}, 
and $\nu\in \nc$.  If $2/(d+3)<\dpp<\min (2/(d-1), 2/3)$, then for some $c>0$,  we have  
\[
\| \sum_j  \sum_{k\sim_\nu k'} \chi_{k}^{\nu} \fP_\lambda[\psi_j] \tilde \chi_{k'}^{\nu}  \|_{p\to p'} \lesssim   \bpqp  \mom^{c}.
\] 
\end{prop}

  Combining  Proposition \ref{lala}  and   \ref{l0}, we prove Proposition \ref{improvedbound}. 
  
  \begin{proof}[Proof of Proposition \ref{improvedbound}] Proposition \ref{lala}  and   \ref{l0}  give  \eqref{eq:la}, from which  \eqref{scaled-c} follows for $2/(d+3)<\dpp<\min (2/(d-1), 2/3)$. 
 Interpolating the estimate  with
 \[\|   \chi_{\mu}^+ \mathfrak P_\lambda  \chi_{\tilde\mu}^+   \|_{2\to 2} \lesssim   \lambda^{-d/2}(\mu\tmu)^{1/4}, \]
 which follows from \eqref{stronger1} after scaling, 
 we obtain \eqref{scaled-c} for $0<\dpp<\min (2/(d-1), 2/3)$.
 \end{proof}

Since $\nu\in \nc$, the estimates in \eqref{b-case-kernel} remain valid. So,  we have   the first and third estimates  in \eqref{eq:sec8}. By the same argument as before we see 
\begin{align*}
\big( \sum_{2^{-j}\ll \sqrt{\mu}} +  \sum_{2^{-j}\gg \sqrt{\mu}}\big) \|  \sum_{k\sim_\nu k'} \chi_{k}^{\nu} \mathfrak P_\lambda [ \psi_j]  \tilde\chi_{k'}^{\nu}\|_{p\to p'}\lesssim \bpqp\mom^c 
\end{align*}
for some $c>0$ because ${2}/{(d+3)}<\dpp< {2}/({d-1})$. 
Thus, by Lemma \ref{kkpsum} the proof of Proposition \ref{l0} is now reduced to showing   
\begin{equation}
\label{bbb011}
\|   \chi_{k}^{\nu} \mathfrak P_\lambda [ \psi_j]  \tilde\chi_{k'}^\nu \|_{p\to p'} \lesssim   \bpqp  \mom^{c}, \quad  \  2^{-j}\sim \smu,  \quad  k\sim_\nu k'.
\end{equation}

In what follows we make further reductions to prove  \eqref{bbb011}.   
By \eqref{rot2},   we may assume that
\begin{align*}
A_{k}^{\nu}\subset\, &\mathfrak R:=\big\{x: |x_1-1|\sim \mu, \,\,|\ol x|\lesssim  \mmt  \big \},
\\
\tilde A_{k'}^{\nu} \subset\, &\tilde {\mathfrak R}:= \big\{y: |y_1-1|\sim {\tilde\mu}, \,\,  |\ol y|\lesssim  \mmt  \big\},
\end{align*}
where  $x=(x_1, \ol x)\in \mathbb R\times \mathbb R^{d-1}$.    For $x^*\in \mathbb R^d$ and $r>0$,  we denote 
\[
\mathfrak R_{\mu, {\tilde\mu}}(x^\ast, r)=\big\{  (x_1, \ol x): |x_1-x_1^\ast|\le r \mu, \,, |\ol x-\ol{x}^\ast|\le r\mmt   \big \}. 
\] 
If $(x,y)\in  \mathfrak R\times \tilde {\mathfrak R}$,  unlike the previous cases,  $\cD(x,y)$ may  vanish.  
To handle this, we  further    localize the value of $\mathcal D$ by decomposing   $\mathfrak R\times \tilde {\mathfrak R}$.  
To this end,  we break  $\mathfrak R\times \tilde {\mathfrak R}$ into finitely many disjoint rectangles  $\mathfrak R_{\mu, {\tilde\mu}}(x',\epsilon)\times \mathfrak R_{{\tilde\mu}, \mu}(y', \epsilon)$ 
so that  $ |\mathcal D|\ll  \vepc\mu{\tilde\mu}$, or $ |\mathcal D|\gtrsim \vepc\mu{\tilde\mu}$ holds  on each of those rectangles for a small $\vepc>0$.  
This particular form of decomposition shall be important  later.

\begin{lem}\label{dychotomy}  Let $\mu, \wmu$ satisfy \eqref{mucon} and $0<\epsilon\le \varepsilon_\circ$. 
Let
$(x^\ast, y^\ast) \in \mathfrak R\times \tilde{\mathfrak R}$. If  $x' \in \mathfrak R_{\mu, {\tilde\mu}}(x^\ast, \epsilon)$ and $y'\in \mathfrak R_{{\tilde\mu},\mu}(y^\ast, \epsilon)$, then   
for a constant $C$ we have  
\Be
\label{det-control} |\mathcal D(x',y')-\mathcal D(x^\ast,y^\ast)|\le  C \epsilon \mm.
\Ee
\end{lem}

\begin{proof} 
Denote $z':=(z_1', \dots, z_{2d}')=(x',y')$ and $z^*:=(z_1^*, \dots, z_{2d}^*)=(x^\ast, y^\ast)$, then set 
$\mathcal D_0=\mathcal D(z')$ and $ \mathcal D_k=\mathcal D(z_1^*,\dots, z_k^*, z'_{k+1}, \dots, z'_{2d}),$ $1\le k\le 2d.$
So, $
 \D(x',y')-\D(x^\ast,y^\ast)=\sum_{k=0}^{2d-1} (\mathcal D_k-\mathcal D_{k+1}).
$
Thus, \eqref{det-control} follows  if we show
\[  | \mathcal D_k-\mathcal D_{k+1}|\le C\epsilon {\tilde\mu}\mu, \quad   1\le k\le 2d. \]
Note  $|x_1'-x_1^\ast|\le \epsilon \mu$, $|y_1'-y_1^\ast|\le \epsilon {\tilde\mu}$, 
and $|x_j'-x_j^*|, |y_j'-y_j^*|\le \epsilon \mmt$, $j\ge 2$. By the mean value theorem  we only need to show
\[  |\partial_{x_1} \cD| \lesssim  \tilde\mu,  \quad  |\partial_{y_1} \cD| \lesssim \mu,  \quad   |\partial_{\ol x} \cD|,   |\partial_{\ol y} \cD| \lesssim \mmt \]    
on $\mathfrak R\times \tilde{\mathfrak R}$. Since $\partial_x\cD (x,y) = 2\inp xy y- 2x$ and $\partial_y\cD (x,y) = 2\inp xy x- 2y$, 
it is clear that $|\partial_{\ol x} \cD|,   |\partial_{\ol y} \cD| \lesssim \mmt $.  
Writing $\partial_{x_1} \cD(x,y) = 2x_1(y_1^2-1)+2\inp{\ol x}{\ol y}y_1$, we see $|\partial_{x_1} \cD| \lesssim  \tilde\mu$. 
Similarly, we get  $|\partial_{y_1} \cD| \lesssim  \mu$.
\end{proof}

Let $\epsilon=c\vepc$ for a small enough $c>0$. 
Let $\mathfrak R_{\mu, {\tilde\mu}}(x_k, \epsilon/2)$, $1\le k\le K$, and $\mathfrak R_{{\tilde\mu},\mu}(y_l, \epsilon/2),  1\le l\le L$ be almost disjoint rectangles  
which  cover the rectangles $\mathfrak R$, $\tilde{\mathfrak R}$, respectively. 
We denote  
\[ \mathcal B=\mathfrak R_{\mu, {\tilde\mu}}(x_k, \epsilon), \quad \tilde \cB=\mathfrak R_{{\tilde\mu},\mu}(y_l, \epsilon).\]
 Taking $c>0$ small enough, by Lemma \ref{dychotomy} 
we may assume that one of the following holds:
\begin{align}
\label{large-d}
  &|\mathcal D(x,y)|\gtrsim \vepc \mm, \  \  \ \,  (x,y)\in   {\mathcal B}\times \tilde {\mathcal B},  
   \\ 
   \label{small-det}
     \ & |\mathcal D(x,y)|\ll \vepc \mm, \  \  \  (x,y)\in   {\mathcal B}\times \tilde {\mathcal B}. 
  \end{align} 
  
  Let   $\tilde \chi_ {\cB}$, $\tilde \chi_ {\tilde\cB}$ be  smooth functions adapted to the rectangles $\cB, \tilde \cB$, respectively, i.e.,  $\supp  \tilde \chi_ {\cB}\subset \mathcal B$, $\supp  \tilde \chi_ {\tilde \cB}\subset \tilde\cB$, and  
  \Be 
\label{d-bound} 
\partial_{x_1}^{\alpha_1} \partial_{\ol x}^{\ol \alpha}  \tilde \chi_ {\cB}=O( \mu^{-\alpha_1}(\mu\tmu)^{-|\tilde\alpha|/2}), \quad 
\partial_{y_1}^{\alpha_1} \partial_{\ol y}^{\ol \alpha}  \tilde \chi_ {\tilde\cB}=O( \tmu^{-\alpha_1}(\mu\tmu)^{-|\tilde\alpha|/2}). 
\Ee 
For the proof of \eqref{bbb011} it is enough to show   
\begin{equation}
\label{bbb01}
\|  \tilde \chi_{{\mathcal B}} \fP_\lambda[\psi_j] \tilde \chi_{\tilde {\mathcal B}} \|_{p\to p'} \lesssim   \bpqp  \mom^{c}, \quad 2^{-j}\sim \smu
\end{equation} for some $c>0$ while either \eqref{large-d} or \eqref{small-det}  holds.  When \eqref{large-d} holds, one  can show \eqref{bbb01}  in the same manner as in the proof of Proposition \ref{lala}. 
Indeed, since   $2^{-\nu}\sim \mmt$,  $2^{-j}\sim \smu$, and $|\mathcal D|\sim  \mm$ on ${\mathcal B}\times \tilde {\mathcal B},$  
we have the second  estimate in \eqref{eq:sec8}, which gives
\eqref{bbb01} for $1\le p<2$.  

Therefore, to complete the proof of Proposition \ref{l0}  it suffices to show  
the following.

\begin{prop} \label{critical-muj}   Let $d\ge 3$ and  $\mu, \wmu$ satisfy \eqref{mucon}.   
Suppose that 
\eqref{small-det} holds. Then,  if $0<\dpp<  2/3$,  we have 
\eqref{bbb01}  for a constant  $c>0$. 
\end{prop}

\subsection{2nd-order derivative of $\cP$}   
To prove Proposition \ref{critical-muj}, we shall dyadically decompose $\fP_\lambda[\psi_j](x,y)$ away from the zero   
of $\partial_s^2\cP$ (see  \eqref{S3-ptkernel}).   A computation shows
\[
 \partial_s^2 \mathcal P(x,y,s) =-\inp xy\frac{\mathcal R(x,y,\cos s) }{ \sin^3 s},
\]
where
\[  \mathcal R(x,y,\tau) =  \tau^2-{\inp xy}^{-1}(|x|^2+|y|^2)\tau+ 1.\]
  For $(x,y)\in \cB\times \tilde \cB$, the equation $\mathcal R(x,y,\tau)=0$ has two distinct  zeros $\tau^\pm(x,y)$: 
\Be
\label{roots}  \tau^\pm(x,y) =\frac{|x|^2+|y|^2\pm |x+y||x-y|}{2\inp xy}.
\Ee
Since   $\tau^+(x,y)>1$,   $\tau^-(x,y)$ is more relevant for our purpose.

\begin{lem} 
\label{scs}  
Define a function  $S_c:\cB\times \tilde \cB\to (0, \pi/2)$  by 
\begin{equation}
\label{sc-}\cos S_c(x,y)= \tau^-(x,y).
\end{equation}
Then,  $S_c$ is smooth and  $S_c(x,y)\sim  |x-y|^{1/2}$ for  $(x,y)\in \cB\times \tilde \cB.$
 \end{lem}
 
We here record a few identities, which are to be  useful later: 
\begin{align}
\label{ph-dd2}
  \partial_s^2 \mathcal P(x,y,s)& =-\inp xy\frac{(\cos S_c(x,y)- \cos s)( \tau^+(x,y)-\cos s) }{ \sin^3 s},
\\[2pt]
\label{cossin}  \tau^+-\cos\sxy &=    \frac{|x+y||x-y|}{\inp xy}= \frac\sut{\cos S_c}.
\end{align}
From now on, for simplicity  we occasionally omit the arguments $(x,y)$ of $S_c$ and related functions  as long as there is  no ambiguity. 

\begin{proof}[Proof of Lemma \ref{scs}] Let $(x,y)\in \cB\times \tilde \cB$. 
From \eqref{roots}  we note  
\Be 
1-\tau^-(x,y)=2^{-1}|x-y|(|x+y|-|x-y|)/\inp xy.
\label{hohoho}
\Ee So,  $\tau^-(x,y)\in [1-c_2\mu, 1-c_1\mu]$ for some positive constants $c_2>c_1>0$.  Thus, it is 
clear that $S_c$ is smooth on $\cB\times \tilde \cB.$  Since  $1-\cos S_c(x,y)\in [c_1\mu, c_2\mu]$, to see $|x-y|^{1/2}\sim  S_c(x,y)$  it suffices to observe that 
\begin{equation}
\label{sc}
\begin{aligned}
1-\cos S_c(x,y) = \frac{2|x- y|}{|x+y|+|x-y|}\,.
\end{aligned}
\end{equation} 
Note $|x+y|^2-|x-y|^2=4\inp xy$. Thus, \eqref{sc} follows by \eqref{hohoho}. 
\end{proof}

\subsection{Reduction to $L^2$ estimates}  From \eqref{sc}  we note 
\[
 |\cos S_c(x,y)-\cos S_c(x^*,y^*)|\lesssim \vepc  \mu, \quad(x,y),\,(x^*,y^*)\in \cB\times \tilde\cB.
\] 
Consequently,  $S_c(\cB\times \tilde {\mathcal B})$ is contained in an interval of length $C_1\vepc\mu^{1/2}$ for a constant $C_1$.
Using this, we make  a further localization.  

Let $\fc_{\cB}$ denote  the center of the rectangle $\cB\times \tilde {\mathcal B}$, and set 
\[
\psi_{\cB}(s)= \psi_\circ\Big(\frac{s-S_c(\fc_\cB)}{C_1\vepc  {\mu}^{1/2}} \Big),
\] 
where  $\psi_\circ\in C_c^\infty(-4,4)$ such that $ 
\psi_\circ=1$ on $[-2,2]$. 
We decompose  
\begin{align*}
\wt\chi_{\mathcal B}
\fP_\lambda [ \psi_j ]  \wt\chi_{\tilde {\mathcal B}}= \fP^0 + \fP^1: = \wt\chi_{\mathcal B}
\fP_\lambda [ \psi_j\psi_{\cB} ]  \wt\chi_{\tilde {\mathcal B}} + \wt\chi_{\mathcal B}
\fP_\lambda [ \psi_j(1-\psi_{\cB})]  \wt\chi_{\tilde {\mathcal B}}.
\end{align*}
The operator $\fP^1$ is easy to handle. In fact, one can show without difficulty 
\begin{align*}
\|\fP^1\|_{2\to 2}\lesssim \lambda^{-\frac d2}(\mu\tilde\mu)^{\frac14}, \qquad
\|\fP^1\|_{1\to \infty}\lesssim \lambda^{-\frac12}\mu^{-\frac{d+1}{4}}.
\end{align*}
Interpolation gives $\| \fP^1\|_{p\to p'} \lesssim \bpqp  \mom^{c}$  for $1\le p<2$.
The $L^2$ bound follows from Lemma \ref{l2-mmj23}. For the  $L^1$--$L^\infty$ bound,  we note  $|\cos s - \cos S_c(x,y)|\sim\mu$ if $s\in\supp(\psi_j(1-{\psi_{\cB}}))$ and $(x,y)\in \cB\times\tilde{\cB}$.
We also have $\tau^{+}(x,y)-\cos s\gtrsim \mu$ for $(x,y)\in \cB\times \tilde \cB$  since $\tau^+(x,y)-1\sim |x-y|$.  Thus, recalling \eqref{ph-dd2}, 
we see $ |\partial_s^2\cP(x,y,s)|\gtrsim  \mu^{1/2}$ for $s\in\supp  (\psi_j(1-{\psi_{\cB}}))$ if  $(x,y)\in \cB\times \tilde \cB$.  Applying the van der Corput lemma to  the integral  $\fP^1(x,y)$ (see \eqref{S3-ptkernel}) gives the desired estimate $|\fP^1(x,y)|\lesssim  \lambda^{-1/2}\mu^{-(d+1)/{4}}$.

\begin{prop}
\label{mm2}  Let $d\ge 3$ and $\mu, \wmu$ satisfy \eqref{mucon}.   If  $2/(d+1)<\dpp<  2/3$, then  for some $c>0$ we have
\begin{equation}
\label{m0}
\|\mathfrak P^0\|_{p\to p'} \lesssim   \bpqp \mom^c. 
\end{equation}
\end{prop}

Now, the proof of Proposition \ref{critical-muj} is  straightforward. 

\begin{proof}[Proof of Proposition \ref{critical-muj}]
Combining the estimates for $\mathfrak P^0$ and $\mathfrak P^1$, we obtain \eqref{bbb01}  for $2/(d+1)<\dpp<  2/3$. Meanwhile,  we have the estimate 
$\|\tilde\chi_{{\mathcal B}} \mathfrak P_\lambda [ \psi_j ]  \tilde\chi_{\tilde {\mathcal B}}\|_{2\to 2}  \lesssim\lambda^{-\frac d2}(\tmu\mu)^{1/4}$ by Lemma \ref{l2-mmj23}. 
Interpolation yields \eqref{bbb01} for $0<\dpp<  2/3$. 
\end{proof}

If $(x,y)\in \cB\times \tilde \cB$,  $\cD(x,y)$ is no longer bounded away from the zero. To get the correct  order of decay  in $\lambda$, i.e.,  
$O(\lambda^{-1/2})$, we consider $\partial_s^2\cP$,  which alone is not enough to give a favorable lower bound since it  also vanishes at some point. 
Such difficulty  is typically circumvented by considering  $\partial_s \cP$ and $\partial_s^2 \cP$ together. 
However,  this is not viable in our situation since the zeros of $\partial_s \cP$ and $\partial_s^2 \cP$ merge as $\cD(x,y)\to 0$. This leads us to  break  the integral  away from $S_c$. 

 Inserting  the cutoff functions  $\tilde\psi(2^{l}(s-S_c))$, we decompose 
 \begin{equation}
\label{mmm}
 \mathfrak P^0=  \ssum{l} \,\mathfrak P_{l}:=  \ssum{l}  \, \tilde\chi_{\cB}\fP_\lambda[\psi_j \psi_{\cB}\tilde\psi(2^{l}(\cdot-S_c)) ] \tilde\chi_{\tilde \cB},
 \end{equation}
 where $\tilde\psi=\psi(|\cdot|)$.  Clearly, we have 
\Be 
\label{interval}
2^{-l}\le \vepc \mu^{1/2},
\Ee
since $\mathfrak P_{l}=0$ otherwise.  As to be  seen later, \eqref{interval} makes it possible to render the minor contribution  manageable  if we take $\vepc$ small enough. 
We set  
\[ \varepsilon_1=C\vepc^{1/2} \]
for a large constant $C>0$.\!  

We have different estimates for  $\mathfrak P_l$ depending on $l$.  

 \begin{lem}\label{propp1}  Let  $d\ge 2$,  and $\mu, \wmu$ satisfy \eqref{mucon}. 
If  $2^{-l}<\varepsilon_1{\wmu}^{1/2}$, then 
\Be\label{p1}
\|\mathfrak P_l\|_{2\to 2}\lesssim \lambda^{-\frac d2} 2^{-l} (\wmu/\mu)^{-\frac12}.
\Ee
 If  $\varepsilon_1{\wmu}^{1/2}\le 2^{-l}$,  then we have
\Be\label{p2}
\|\mathfrak P_l\|_{2\to 2}\lesssim \lambda^{-\frac d2} 2^{\frac l2} \mu^{\frac14}\wmu^{\frac12}.
\Ee
\end{lem}

We postpone the proof of Lemma \ref{propp1} until the next section.  Assuming this  for the moment,
 we prove Proposition \ref{mm2}.    To do this, we set  \[S_c^l(x,y,s) = 2^{-l}s + S_c(x,y).\]  
Then, changing variables $s\to S_c^l(x,y,s)$, we  have
\Be 
\label{kernel-pl}
\mathfrak P_{l}(x,y)=   \chi_{\cB}(x)\chi_{\tilde \cB}(y)  2^{\frac d2j} 2^{-l}  \int    (2^{-\frac d2j}\mathbf a \psi_j \psi_\cB)(\sxyl )\,  \tilde \psi(s)\,  e^{i\lambda  \mathcal P(x,y,S_c^l ) } ds. 
\!\!\!
\Ee 
 Here,  we also drop the arguments of $S_c^l$  for simplicity as before. 
 \begin{proof}[Proof of Proposition \ref{mm2}] 
 Since $|\cos S_c^l-\cos S_c|\sim\mu^{1/2}2^{-l}$ on the support of $ (\mathbf a \psi_j \psi_\cB)\circ \sxyl\,\tilde \psi $, by \eqref{ph-dd2}   we have 
 $|\partial_s^2 \big(\mathcal P(x,y,S_c^l)\big)|\sim 2^{-3l}$. Applying the van der Corput lemma,  for $l$ satisfying \eqref{interval} we get 
 \Be\label{pl1inf}
\|\mathfrak P_l\|_{1\to\infty} \lesssim \lambda^{-\frac12} 2^{\frac l2}  \mu^{-\frac d4}.
\Ee 
  
Interpolation  with \eqref{p2} gives 
\begin{align*}
\|\mathfrak P_l\|_{p\to p'}\lesssim \lambda^{\frac{d-1}{2}\delta(p,p')-\frac d2}2^{\frac l2}\mu^{\frac14-\frac{d+1}{4}\delta(p,p')}\wmu^{\frac12(1-\delta(p,p'))}
\end{align*}
when  $\varepsilon_1{\wmu}^{1/2}\le 2^{-l} \le \vepc {\mu}^{1/2}$. Hence, by summation over $l$ we have
\[
\| \ssum{\varepsilon_1{\wmu}^{1/2}\le 2^{-l} \le \vepc {\mu}^{1/2}}\mathfrak P_l \|_{p\to p'}\lesssim \mathrm B_{p}(\mu,\wmu)({\wmu}/{\mu})^{\frac{d-1}{8}\delta(p,p')}.
\]
By interpolation between the estimates  \eqref{pl1inf} and \eqref{p1} we get
\[
\|\mathfrak P_l\|_{p\to p'}\lesssim \lambda^{\frac{d-1}{2}\delta(p,p')-\frac d2} 2^{-l(1-\frac 32\delta(p,p'))} (\mu\wmu)^{-\frac d8\delta(p,p')} ({\wmu}/{\mu})^{-\frac12+\frac{d+4}{8}\delta(p,p')}
\]
for $2^{-l}<\varepsilon_1\wmu^{1/2}$, which yields
\[
\|\ssum{ 2^{-l}< \varepsilon_1{\wmu}^{1/2}}\mathfrak P_l\|_{p\to p'}\lesssim  \mathrm B_{p}(\mu,\wmu)({\wmu}/{\mu})^{-\frac14+\frac{d+1}{8}\delta(p,p')}
\]
if $\delta(p,p')<{2}/{3}$.  Recalling  \eqref{mmm}, we combine the estimates above and  obtain  \eqref{m0}
for $2/(d+1)<\delta(p,p')<2/3$.
\end{proof}

\section{$L^2$ estimate: Proof of Lemma \ref{propp1}}
\label{sec:l2}
 For $a\in C_c^\infty(\mathbb R^d\times \mathbb R^d)$ and a smooth function $\phi$ on  $
\supp a$,  
we define 
\[ \mathcal O_\lambda[\phi, a]f(x)=\int e^{i\lambda  \phi(x,y) } a(x,y) f(y) dy.\]
We denote $\partial_x f=( \partial_{x_1} f, \dots, \partial_{x_d}f)^\intercal$ and  $\partial_x^\intercal f= ( \partial_{x_1} f, \dots,   \partial_{x_d} f)$, so that $\partial_x \partial_y^\intercal=(\partial_{x_i}\partial_{y_j})_{1\le i,j\le d}$. 
  We now recall the following lemma.

\begin{lem}[{\cite{hor}, \cite[p.377]{St93}}]
\label{generalized}  
Let $\lambda>0$.  If $ \det( \partial_x \partial_y^\intercal \phi )\neq 0$ on $\supp a$, then there is a constant $C=C(\phi,a)$ such that 
\[\|\mathcal O_\lambda[\phi, a]f\|_2\le C\lambda^{-\frac d2}\|f\|_2.\]
\end{lem}

From \eqref{kernel-pl},  we see 
\Be\label{frakpl}
\mathfrak P_l f= 2^{\frac d2j} 2^{-l} 
\int  \tilde \psi(s)  \mathcal O_\lambda  [\Phi_s, A_s] f\,ds,
\Ee
where
\begin{align}
    \Phi_s(x,y) &= \mathcal P(x,y,S_c^l), \label{Ps}\\
    A_s(x,y)&= \tilde \chi_ {\cB}(x) \wt\chi_{\tilde \cB}(y)   \big ( 2^{-\frac d2j} \mathfrak a\psi_j \psi_\cB \big)(\sxyl).\label{As}
    \end{align}
    To obtain the estimates \eqref{p1} and \eqref{p2}, one may try to use Lemma \ref{generalized} for $\mathcal O_\lambda  [\Phi_s, A_s]$.
However,   $\mathcal O_\lambda [\Phi_s, A_s]$ exhibits different natures  depending on $\mu,\tmu$ and $l$.  
When $\mu\sim \tmu$,  via  a  suitable change of variables $\mathcal O_\lambda [\Phi_s, A_s]$  can be handled by the estimate in  Lemma \ref{generalized}. 
However, when ${\tilde\mu}/\mu$ gets smaller, Lemma \ref{generalized} is not enough (see Lemma \ref{generalized-small2} below).

We  set $\mathbb J=\{ s: 1/4\le |s|\le 1\}$.   
Since $\supp \tilde \psi\subset \mathbb J$ and $2^{-j}\sim \sqrt\mu$, by \eqref{frakpl} the estimate \eqref{p1} follows  from the next lemma.

\begin{lem}\label{key-osc1} Let $d\ge 2$ and  \eqref{mucon} hold.   
If $2^{-l} <\varepsilon_1\wmu^{1/2}$, then    
\Be \label{key-l21}
\|\mathcal O_\lambda [\Phi_s, A_s]\|_{2\to 2}\le C   \lambda^{-\frac d2}  \mu^\frac d4 \mom^{-\frac12 }, \quad \ \  s\in \mathbb J.
\Ee
\end{lem}

When  $\varepsilon_1\wmu^{1/2} \le 2^{-l} \le \vepc \mu^{1/2} $, in order to prove \eqref{p2}
we use  a different expression of $\mathfrak P_l$ so that we can exploit a lower bound on $\partial_s\Phi_s$.
Recalling \eqref{frakpl}, by   integration by parts in $s$  we have  
\[
\mathfrak P_l f(x) = 2^{\frac d2j}2^{-l}\int\int e^{i\lambda\Phi_s(x,y)} \tilde  A_s (x,y) 
ds\, f(y)dy,
\]
where 
\[ \tilde  A_s (x,y) = \wt\psi(s) A_s(x,y) \frac{\partial^2_s\Phi_s(x,y)}{i\lambda(\partial_s\Phi_s(x,y))^2}-\frac{\partial_s(\wt\psi(s)A_s(x,y))}{i\lambda\partial_s\Phi_s(x,y)}. \]
Note   $\partial_s\Phi_s\neq 0$ if  $A_s\neq 0$. In fact,  we have 
\Be 
\label{lower-}
|\partial_s\Phi_s(x,y)|\gtrsim 2^{-3l},  \quad   (s,x,y) \in \mathbb J\times\cB\times {\tilde \cB}.
\Ee
Using \eqref{sc-} and \eqref{roots},  we see
 $|\inp xy-\cos S_c(x,y)| \lesssim |\mathcal D(x,y)|/|x-y|$. Thus, by \eqref{small-det} it follows  that 
 \Be 
 \label{haha}
 |\inp xy-\cos S_c(x,y)|\le  \vepc\tilde\mu, \quad   (x,y) \in \cB\times {\tilde \cB}.
 \Ee
Note $|\cos S_c^l - \cos S_c|\sim 2^{-l} \mu^{1/2}$, so $|\inp xy-\cos \sxyl(x,y)|\gtrsim 2^{-l}\mu^{1/2} \ge \varepsilon_1\mmt.$ 
By \eqref{q-def}, \eqref{small-det}, and our choice of $\varepsilon_1$, we get  $\mathcal{Q}(x,y,\cos \sxyl)  \gtrsim 2^{-2l}\mu$. 
So, \eqref{lower-} follows  by  \eqref{ph-d} since $\partial_s \Phi_s=2^{-l}\partial_s \mathcal P(x,y,S_c^l).$

\medskip

The estimate \eqref{p2} is an immediate consequence of the following.

\begin{lem}\label{key-osc2} Let $d\ge 2$ and  \eqref{mucon} hold. If 
$\varepsilon_1\wmu^{1/2} \le 2^{-l} \le \vepc \mu^{1/2} $, then 
\Be \label{key-l22}
\|\mathcal O_\lambda [\Phi_s, \tilde A_s]\|_{2\to 2}\le C   \lambda^{-\frac d2} \mu^\frac {d+1}4 \wmu^{\frac12}  
 2^{\frac32 l} , \quad \ \  s\in \mathbb J.
\Ee
\end{lem}

 For  the rest of this section, we assume \eqref{mucon} and \eqref{interval} hold, and $(x,y)\in \cB\times {\tilde \cB},$ $s\in \mathbb J$, even if it is not 
 made explicit.

 \subsection{Bounds on $\partial_{x}^\alpha\partial_y^\beta \Phi_s$ and $\partial_{x}^\alpha A_s$} In order to prove Lemma 
\ref{key-osc1} and \ref{key-osc2}, we need  estimates for the derivatives of $\Phi_s$ and $A_s$. 
For a given multi-index $\alpha\in \mathbb N_0^{d}$, we write $\alpha=(\alpha_1, \ol \alpha)\in \mathbb N_0\times \mathbb N_0^{d-1}.$

\begin{lem}  Let $s\in \mathbb J$ and $(x,y) \in \cB\times {\tilde \cB}$. Suppose  \eqref{mucon}  holds. 
 Then,
\begin{align}
\label{est-partial}
|\partial_{x}^\alpha\partial_y^\beta S_c(x,y)|&\lesssim  \mu^{\frac12-|\alpha|-|\beta|},  
\\
\label{est-amp}
|\partial_{x}^\alpha A_s(x,y)|&\lesssim  \mu^{-\alpha_1}(\mm)^{-\frac{|\ol \alpha|}{2}}  .
\end{align} 
\end{lem}

\begin{proof} To show \eqref{est-partial}, we note from \eqref{sc} that 
$S_c(x,y)$ behaves like $g(x,y):=|x-y|^{1/2}$, and   that  \eqref{est-partial} and \eqref{est-amp} are easy to show if $S_c$ is replaced by $g$.

We prove \eqref{est-partial} in an inductive way by making use of \eqref{sc}.  By Lemma \ref{scs},    \eqref{est-partial} holds with $\alpha=\beta=0$  since $\dist(\cB, {\tilde\cB})\sim \mu$. 
We now assume  \eqref{est-partial} is true for $|\alpha|+|\beta|\le N$. Applying $\partial_{x}^\alpha\partial_y^\beta $ on both side of \eqref{sc} for $|\alpha|+|\beta|=N+1$,  it is not difficult to see 
\Be
\label{derivative2}
\sin \sxy(x,y) \partial_{x}^\alpha\partial_y^\beta S_c(x,y) +O\big(\mu^{-N}\big)= O\big(\mu^{-N}\big)
\Ee
for $(x,y)\in \cB\times \tilde\cB.$ Indeed, the right hand side of \eqref{sc} behaves  as if  it were $|x-y|$. On the  left hand side of  \eqref{derivative2}, the terms other than the first one are given by a linear combination of 
the products of   $\sin \sxy$, $\cos \sxy$, and  $\prod_{i=1}^l\partial_{x}^{a_i}\partial_y^{b_i} S_c$ with $l\ge 2$ and $\sum_{i=1}^l ( |a_i|+|b_i|)\le N+1$. Our induction assumption 
shows those  are $O(\mu^{-N})$, thus we get \eqref{derivative2}. 
Since  $|x-y|\sim \mu$,  $\sxy(x,y)\sim \smu$ by Lemma \ref{scs}. Therefore,  \eqref{derivative2} gives  $\partial_{x}^\alpha\partial_y^\beta S_c(x,y)=O(\mu^{-N-1/2})$ as desired. 

Once we have 
\eqref{est-partial}, \eqref{est-amp} is easier to show. 
Recall \eqref{As}. Since we have  \eqref{d-bound},  it is sufficient to show  
\[ \partial_{x}^\alpha\big( ( 2^{-jd/2} \mathfrak a\psi_j\psi_{\cB})\circ \sxyl \big)=O( \mu^{-|\alpha|}).\]
Note $(d/ds)^k(2^{-jd/2} \mathfrak a\psi_j)=O(\mu^{-k/2}) $ and  $(d/ds)^k  \psi_{\cB}=O(\mu^{- k/2})$. By  the chain rule and \eqref{est-partial},  
we get $\partial_{x}^\alpha ( (2^{- jd/2} \mathfrak a\psi_j)\circ \sxyl )=O(\mu^{-|\alpha|})$ and $\partial_{x}^\alpha (\psi_{\cB}\circ \sxyl )= O( \mu^{-|\alpha|})$.
From those estimates the desired bound follows. 
\end{proof}

\begin{lem} \label{lem:ph-derivative}  Let $s\in \mathbb J$ and $(x,y) \in \cB\times {\tilde \cB}$. Suppose that   \eqref{mucon} and \eqref{interval} hold. 
  Then,  
\begin{equation}
\label{ph-derivative}
 \partial_x^\alpha \partial_y^\beta  \Phi_s(x,y)=O(\mu^{\frac32-|\alpha|-|\beta|}),  \quad \ |\alpha|, |\beta|\ge 1. 
 \end{equation}
\end{lem}

\begin{proof} Let us set
\[ u(t)=\frac{t-\sin t}{2}, \qquad \,  v(t)= \frac{2\cos t -2+\sin^2 t}{2\sin t}, \qquad \,  w(t)=\frac{\cos t}{2\sin t} .\]
Recalling \eqref{Ps}  and \eqref{p-def}, we  write 
 \begin{align*}
 \Phi_s(x,y)= u\big(\sxyl\big)  +  2^{-1}(1-\inp xy) \sin\sxyl + \inp xy \,v\big(\sxyl \big) + {|x-y|^2} \,w\big(\sxyl\big) .
\end{align*}
From \eqref{est-partial} it is clear that  $|\partial_{x}^\alpha\partial_y^\beta \sxyl (x,y)|\lesssim  \mu^{\frac12-|\alpha|-|\beta|}$. We note 
$(d/dt)^k u(t),$ $(d/dt)^k v(t)=O(t^{3-k})$, and $(d/dt)^k w(t)=O(t^{-1-k})$. Since  $\sxyl(x,y)\sim \smu$, $1-\inp xy=O(\mu)$,  and $|x-y|\sim \mu$,  a routine computation 
yields \eqref{ph-derivative}. 
\end{proof}

\begin{rem}\label{derivs} Under the same assumption as in Lemma \ref{lem:ph-derivative},  
we also have 
$|\partial_{x}^\alpha\partial_y^\beta \sin^m S_c| ,$  $  |\partial_{x}^\alpha\partial_y^\beta \sin^m \sxyl | \lesssim  \mu^{\frac m2-|\alpha|-|\beta|}$ for  any $m\in \mathbb Z$, 
and $|\partial_{x}^\alpha\partial_y^\beta (1-\cos S_c)|$, $|\partial_{x}^\alpha\partial_y^\beta (1-\cos\sxyl)|\lesssim   \mu^{1-|\alpha|-|\beta|}$ for any multi-indices $\alpha, \beta$. 
\end{rem}

\subsection{Estimate for $\partial_x\partial_y^\intercal \Phi_s$}
For a given matrix $\mathbf N$, we denote by $\mathbf N_{i,j}$ the $(i,j)$-th entry of $\mathbf N$. 
We consider a $d\times d$ matrix   $\mathfrak M^0$  which is given   as follows: 
\begin{align*}
\mathfrak M^0_{i,i}(x,y)&= 1-\frac{(x_i-y_i)^2}{|x-y|^2}, &&  1\le i \le d,
\\
 \mathfrak M_{1,j}^0 (x,y)&= -\frac{x_j-y_j}{2|x-y|^2}\big(2x_1-(1+\cos S_c)y_1\big), &&   j\ge 2,  
\\[2pt]
 \mathfrak M_{i,1}^0(x,y)&= \mathfrak M_{1,i}^0(y,x), &&   i\ge 2, 
\\[2pt]  
\mathfrak M_{i,j}^0(x,y) &= -\frac{(x_i-y_i)(x_j-y_j)}{|x-y|^2}, &&  i,j\ge 2,\,\,  i\neq j . 
\end{align*} 

 The following lemma  shows the matrix $\partial_x \partial_y^{\intercal}   \Phi_s(x,y)$ is  close to $\mathfrak M^0$.  
 Let   \[ 
E^\alpha_k =\begin{cases}
   \, \vepc(\wmu/\mu)^\frac1k + 2^{-2l}\mu^{-1}, & \, \alpha = 0, 
    \\[3pt]
 \,   \big((\wmu/\mu)^\frac1k + 2^{-2l}\mu^{-1}\big) \mu^{-\alpha_1} (\mu\wmu)^{-\frac{|\ol\alpha|}{2}},  & \, \alpha\neq 0, 
    \end{cases} \qquad k=1,2.
\]  
 
 \begin{lem}\label{matlem}  
Let $s\in \mathbb J$ and $(x,y) \in \cB\times {\tilde \cB}$.  Set
\[
\mathfrak M(x,y)=-  \sin S_c^l(x,y) \partial_x \partial_y^{\intercal}   \Phi_s(x,y).
\] Suppose that   \eqref{mucon} and \eqref{interval} hold. 
Then, $\mathfrak M=\mathfrak M^0+\mathcal E$ and 
$\CE_{i,j}$ satisfies  
  \begin{align}
  \label{erbd}
  | \partial_x^\alpha \CE_{i,j}(x,y) | \lesssim 
    \begin{cases} \  E^\alpha_1,  \qquad  (i,j)=(1,1)\text{ or }i,j\ge 2,
       \\[2pt]
        \   E^\alpha_2,  \qquad  \qquad\quad  \text{otherwise}. 
           \end{cases}
              \end{align}  
            
\end{lem}

Note $E_1^\alpha\le E_2^\alpha$. 
We postpone the proof of Lemma \ref{matlem} until the end of this section.
Instead, we deduce a couple of lemmas from it 
for later use.  By $\tilde {\mathfrak M}$ we denote the $(1,1)$ minor matrix of $ \mathfrak M$, i.e.,  $\tilde {\mathfrak M}=(\mathfrak M_{i+1,j+1})_{1\le i,j\le d-1}$. 

\begin{lem}\label{matdet:lem}
Let $s\in \mathbb J$ and $(x,y) \in \cB\times {\tilde \cB}$. Suppose   \eqref{mucon} and \eqref{interval} hold. Then, 
we have 
\begin{align}
\label{matdet}
\det{\mathfrak M}(x,y)&\sim {\wmu}/{\mu},  \qquad\ \ 2^{-l}\le\varepsilon_1\wmu^{1/2},
\\
\label{matdet2}
\det  \tilde{\mathfrak M}(x,y)  &\sim 1, \qquad \qquad 2^{-l}\le\vepc\mu^{1/2}.
\end{align}
\end{lem}

\begin{proof}  
Note $\mathfrak M_{i,i}^0 = 1 + O(\wmu/\mu)$ for $2\le i\le d$ and  $\mathfrak M_{i,j}^0= O(\tilde{\mu}/\mu)$ for $2\le i\neq j\le d$. 
Thus,  by Lemma \ref{matlem}  we have $ \tilde {\mathfrak M}=\mathbf I_{d-1} + O(\vepc+ \tilde{\mu}/\mu)$. Since $\tilde{\mu}/\mu\le \vepc$,  \eqref{matdet2} follows if 
$\vepc$ is small enough. 

For the proof of \eqref{matdet},  we may assume 
\begin{align}\label{coordinates}
    x=(r,0,0,\dots, 0),\quad y=(\rho,h,0,\dots, 0). 
\end{align} 
Indeed, observe  $S_c(x, y)=S_c(\mathbf U x, \mathbf Uy)$ and 
$\Phi_s(x, y)=\Phi_s(\mathbf Ux, \mathbf Uy)$ for $\mathbf U\in \mathrm O(d)$. The second identity and a computation  show
$\partial_x \partial_y^\intercal \Phi_s(x, y)$ $= \mathbf U^\intercal \partial_x \partial_y^\intercal \Phi_s(\mathbf Ux, \mathbf Uy)  \mathbf U,$ 
i.e., $\mathbf U  \mathfrak M(x,y) \mathbf U^\intercal=\mathfrak M(\mathbf Ux, \mathbf Uy)$. 
Since $\det{\mathfrak M(\mathbf Ux, \mathbf Uy)}$ $= \det{\mathfrak M(x,y)}$,
we need only to choose $\mathbf U\in \mathrm O(d)$ such that $\mathbf Ux=(r,0,0,\dots, 0)$, $\mathbf Uy=(\rho,h,0,\dots, 0)$.

Let $\mathbf M$ be the $2\times 2$ matrix  given by 
\begin{align*}
\mathbf M=\frac1{2|x-y|^2} \begin{pmatrix} 2h^2 & \  (1+\cos S_c)\rho h-2rh 
 \\
2\rho h -(1+\cos S_c)rh  &  \   2(r-\rho)^2
\end{pmatrix} .  
\end{align*}
By Lemma \ref{matlem} and  \eqref{coordinates},  the matrix $\mathfrak M$ is of the form 
\begin{align*}
    \mathfrak M (x, y)=  \begin{pmatrix}
    \mathbf M && 0 \\
    0 && \mathbf I_{d-2}
    \end{pmatrix} + \tilde \CE, 
     \end{align*}
   where $\tilde \CE_{i,j} =O((\vepc+\varepsilon_1)(\wmu/\mu) )$  if $(i,j)=(1,1)\text{ or }i,j\ge 2$,   and  $\tilde \CE_{i,j} =O((\vepc+\varepsilon_1)(\wmu/\mu)^{1/2} )$ otherwise. 
Since $(x,y)\in \cB\times \tilde\cB$,  we have  $1-r\sim \mu$, $1-\rho\sim \tilde\mu$, and $h\sim (\mu\tilde\mu)^{1/2}$.  Hence,  $\mathbf M_{1,1}\sim \wmu/\mu$,   $\mathbf M_{2,2}\sim  1,$  and   $\mathbf M_{1,2}$, $\mathbf M_{1,2}=O((\wmu/\mu)^{1/2}  ) $. Consequently,  we see  $\det \mathfrak M(x, y)= \det \mathbf M+ O((\vepc+\varepsilon_1)\wmu/\mu)$.  
Note $r-\rho\sim \mu$ and $1-\cos S_c \sim \mu$. So, we have
\begin{align*}
\det{\mathbf M} &= \frac{h^2}{4|x-y|^4}(1-\cos S_c)\big(2(r-\rho)^2+r\rho (1-\cos S_c)\big) \sim \frac{\wmu}{\mu}. 
\end{align*}
Therefore, \eqref{matdet} follows  if $\vepc$ is sufficiently small. 
\end{proof}

The following is a straightforward consequence of  Lemma \ref{matlem}. 

\begin{lem}\label{lembd}
Let $s\in \mathbb J$ and $(x,y) \in \cB\times {\tilde \cB}$.  
Suppose  \eqref{mucon} holds.  If $2^{-l}\le\varepsilon_1\wmu^{1/2}$, then  for  $\alpha\in \N_0^d$ 
\begin{align*}
 |\partial_x^\alpha\mathfrak M_{i,j}(x,y)|\lesssim 
   \begin{cases}    
   (\wmu/\mu)\mu^{-\alpha_1}(\mu\wmu)^{-\frac{|\ol\alpha|}{2}}, \quad  & i=j=1, 
   \\[2pt]
     (\wmu/\mu)^\frac12\mu^{-\alpha_1}(\mu\wmu)^{-\frac{|\ol\alpha|}{2}},   \quad  & \text{$i$ or $j = 1$}, 
    \\[2pt]
    \mu^{-\alpha_1}(\mu\wmu)^{-\frac{|\ol\alpha|}{2}},  \quad  &i,j\ge 2.
    \end{cases} 
\end{align*}
Furthermore, the last bound remains valid even if $2^{-l}\le\vepc\sqrt{\mu}$.
\end{lem}

\subsection{Estimate for $\mathcal O_\lambda[\Phi_s,A_s]$ when $2^{-l}< \varepsilon_1\wmu^{1/2}$}
In order to prove the estimate \eqref{key-l21}, we use the following  lemma, which 
differs from  the typical one (see Lemma \ref{generalized}),  in that 
the phase and amplitude functions depend on a parameter. 
We denote $\mathbb B_d(x,r)=\{y\in \mathbb R^d: |y-x|<r\}. $ 

\begin{lem}
\label{generalized-small2}  Let $0<\omega\le 1$.  Let $a$ be a smooth function supported in $S_\omega:= \mathbb B_{d}(0,1)\times  \mathbb B_1(0, \omega)\times  \mathbb B_{d-1}(0,1)$ such that  $
 |\partial_{x}^{\alpha}  a|\lesssim 1$ for $|\alpha|\le d+1$. 
Suppose that  
$|\!\det \partial_x\partial_y^\intercal \phi|\sim 1$ 
 and 
\Be
\label{p-con}
| \partial_{x}^{\alpha}\partial_{y}^{\beta}  \phi| \lesssim  
\begin{cases} 1, &  |\beta| = 1 \\
\omega^{-1 +\frac{|\ol \alpha|+|\ol \beta|}2}, &  |\beta|= 2
\end{cases}
\Ee 
for $1\le |\alpha|\le d+1$ on $S_\omega$. Then, 
we have 
\Be \label{l2-oa} 
\|\mathcal O_\lambda [\phi, a]  f\|_2\lesssim \lambda^{-\frac d2}\|f\|_2. 
\Ee 
\end{lem}

\begin{proof}[Proof of Lemma \ref{generalized-small2}]  
By finite decomposition and translation we may replace $S_\omega$ with 
$\mathbb B_{d}(0,\epsilon_0)\times [-\epsilon_0\omega, \epsilon_0\omega]\times \mathbb B_{d-1}(0,\epsilon_0)$ for a small enough $\epsilon_0>0$. 

We set 
\[ \Psi(x)=  \phi(x,y)-\phi (x, y'), \quad
\quad A(x) =   a(x,y)\ol {a}(x,y'),\] 
and  consider  the integral 
\[ I_\lambda(y,y')= \int e^{i\lambda \Psi(x) }  A(x)dx,
\]
which is the kernel of the operator $\mathcal O_\lambda [\phi, a]^*\mathcal O_\lambda [\phi, a]$.   
The estimate  \eqref{l2-oa} follows  by a standard argument if we  show  
\begin{equation}
\label{iii}
| I_\lambda(y,y')| \le C (1+\lambda|y-y'|)^{-d-1}.
 \end{equation}

Assuming, for the moment, that 
\begin{align}
|\nabla  \Psi(x)|&\gtrsim |y-y'|, 
\label{lowerbdd}
\\
|\partial_x^\alpha   \Psi(x)|&\lesssim |y-y'|,   \quad \forall |\alpha|\ge 2,
 \label{upperbdd}
\end{align}
we prove \eqref{iii}.  Note $\partial^\alpha_x A=O(1)$ for $|\alpha|\le d+1$ since $
 |\partial_{x}^{\alpha}  a|\lesssim 1$. 
Using $\big(\frac{\nabla  \Psi(x)}{i\lambda |\nabla   \Psi(x)|^2}  \cdot \nabla\big) e^{i\lambda \Psi(x) }=e^{i\lambda \Psi(x) }$, by  integration by parts $d+1$ times  we have 
\[ |I_\lambda(y,y')|\lesssim    \sum_{\ell =0}^{d+1} \ \  \sum_{\sum |a_i| \le d+\ell+1, \ |a_i|\ge 2  } \ 
\int \frac{\prod_{i=1}^\ell |\partial_x^{a_i} \Psi(x)|}{\lambda^{d+1}|\nabla  \Psi(x)|^{d+\ell+1}}  \mathcal A(x) dx,    \] 
where $\mathcal A$ is a bounded function supported in $\mathbb B_d(0,1)$.  By  \eqref{lowerbdd}  and \eqref{upperbdd}  
the estimate \eqref{iii} follows.   

We now show \eqref{lowerbdd} and \eqref{upperbdd}.
For \eqref{lowerbdd}
it is sufficient to show \[ |\nabla_x \phi(x,y)- \nabla_x \phi( x,y')|
\gtrsim |y-y'|.\]  By Taylor series expansion, $\nabla_x \phi(x,y)- \nabla_x \phi(x,y')$ is equal to
\[\textstyle \partial_x\partial_y^\intercal \phi (x,y')(y-y')
+\sum_{i, \beta:\, |\beta|=2} O(E_{i,\beta}|(y-y')^\beta|), \] 
where $E_{i,\beta}=\sup_{(x,y)\in S_\omega} |\partial_{x_i}\partial_{y}^\beta  \phi(x,y)|$. Since 
$|y_1-y'_1|\le \epsilon_0 \omega$, by \eqref{p-con}  
 it is clear that $E_{i,\alpha}|(y-y')^\alpha|=O(\epsilon_0|y-y'| )$.  We also have $|\partial_x\partial_y^\intercal \phi (x,y')(y-y')|\sim |y-y'|$ because $(\partial_x\partial_y^\intercal \phi)_{i,j}=O(1)$ and $|\!\det \partial_x\partial_y^\intercal \phi|\sim 1$. Thus,  $ |\nabla_x \phi(x,y)- \nabla_x \phi(x,y')|\gtrsim |y-y'|.$  
For \eqref{upperbdd},  by the mean value theorem  we only need to show $\partial_x^\alpha \partial_{y_j} \phi(x,y)=O(1)$ for $|\alpha|\ge 2$. This follows  from \eqref{p-con}.
\end{proof}

We prove Lemma \ref{key-osc1} by combining 
Lemma  \ref{matdet:lem}, \ref{lembd},  and \ref{generalized-small2}.

 \begin{proof}[Proof of Lemma \ref{key-osc1}]    
 We first transform  $\mathcal O_\lambda [\Phi_s, A_s]$ via scaling and translation so that  we can  apply Lemma \ref{generalized-small2}. 
 Recalling  that $\fc_\cB$ denotes the center of the rectangle $\cB\times \tilde\cB$, 
 we  set 
\begin{align*}
L(x,y)= (L_1x,L_2 y)+  \fc_\cB:=\big( \mu  x_1,\mmt\,\,\ol x, \mu y_1, \mmt\,\, \ol y\big)+ \fc_\cB.\end{align*}
Changing   variables $(x,y)\to L(x,y)$, we  have 
\Be 
\label{op-norm2}
\|\mathcal O_\lambda [\Phi_s, A_s]\|_{2\to 2}  =(\mu/\wmu)^{\frac12}(\mu\wmu)^{\frac d2}
\|\mathcal  O_{\smu\wmu\lambda} \big[ \wt\Phi,  \wt A\,]\|_{2\to 2},
\Ee
where 
\[  
\tilde \Phi(x,y)=(\smu\wmu)^{-1}\Phi_s(L(x,y)) ,\qquad 
 \tilde{A}(x,y)=  A_s(L(x,y)) .
\]
Since $\cB$, $ \tilde \cB$ are of dimensions  about $\vepc\mu\times\vepc\mmt\times\cdots\times\vepc\mmt$,  $\vepc\wmu\times\vepc\mmt\times \cdots \times\vepc\mmt$, respectively,    $\wt A$ is supported in 
$ \mathbb B_{d}(0, 1) \times  \mathbb B_1(0, \wmu/\mu)\times  \mathbb B_{d-1}(0,1).$ 
From \eqref{est-amp} it follows that 
$ \partial_x^\alpha \wt A=O(1)$ for all $\alpha$. 
Since $ \sin S_c(x,y)\sim \smu$,  Lemma \ref{matdet:lem} gives $|\det \partial_x  \partial_y^\intercal   \wt \Phi  |\sim 1$. 
 Indeed, note that 
 \[ \partial_x \partial_y^\intercal \wt{\Phi}= - \frac{\smu}{\sin \sxyl (L(x,y))} (\mu\wmu)^{-1} L^\intercal_1\, \mathfrak M(L(x,y))\,L_2. \]
 Besides,  since $ \partial_x^\alpha (\smu/\sin\sxyl)=(\mu^{-|\alpha|}),$ by Lemma \ref{lem:ph-derivative} and  \ref{lembd} we also have
\begin{align*} 
 \partial_x^\alpha \partial_y^\beta\wt{\Phi} &=\begin{cases}
   O(1), & |\beta| = 1, \\
O\big((\wmu/\mu)^{-1+\frac{|\ol\alpha|+|\ol\beta|}{2}}\big), & |\beta|= 2,
\end{cases}
\end{align*}
 for $(x,y)\in \supp \wt A$ and $\alpha \in \N^d_0$.  In fact, we use  Lemma \ref{lembd} for $|\beta|=1$ and Lemma \ref{lem:ph-derivative} for $|\beta|= 2$. 
Therefore,  we may use Lemma \ref{generalized-small2} with $\omega = \wmu/\mu$ for $ \mathcal  O_{{\smu\wmu\lambda}}[ \wt \Phi,  \wt A\,]$, and we hence  get
\[
\|\mathcal  O_{\smu\wmu\lambda}\big[ \wt\Phi,  \wt A\,\big]\|_{2\to 2}
\lesssim 
\lambda^{-\frac d2}  \wmu^{-\frac d2}\mu^{-\frac d4}.
\] 
By \eqref{op-norm2}   the estimate \eqref{key-l21} follows as desired. 
\end{proof}

\subsection{Estimate for $\mathcal O_\lambda[\Phi_s, \tilde A_s]$ when $2^{-l}\ge \varepsilon_1\wmu^{1/2}$} In this case, as seen in Lemma \ref{matlem}, the matrix 
$\partial_x \partial_y^\intercal \Phi_s$ can be singular. So, we need an approach different  from that used for the case $2^{-l}< \varepsilon_1\wmu^{1/2}$. 
We consider an operator given by 
freezing $x_1, y_1$ (see \emph{Proof of Lemma \ref{key-osc2}} below) and then make use of  \eqref{matdet2}.  

To prove Lemma \ref{key-osc2}, we need the following. 

\begin{lem}\label{bdbs}
Let $2^{-l}\ge \varepsilon_1\wmu^{1/2}$ and $(x,y)\in \cB\times \tilde\cB$. Then, 
\Be
\label{e_bdBs}
|\partial_{\ol x}^{\ol\alpha} \tilde  A_s (x,y)|\lesssim (\lambda 2^{-3l})^{-1}(\mu\wmu)^{-{|\ol\alpha|}/{2}},  \qquad \ol \alpha\in \N_0^{d-1}.
\Ee
\end{lem}

Assuming this for the moment,  we prove Lemma \ref{key-osc2}. 

\begin{proof}[Proof of Lemma \ref{key-osc2}] 
Fixing $s\in \mathbb J$, we set 
\[ 
\bar\Phi (x,y)=(\smu\wmu)^{-1}\Phi_s(\bar L(x,y)) , \qquad 
 \bar{B}(x,y)=  \lambda 2^{-3l}\tilde  A_s (\bar L(x,y)), 
 \]
where 
\[ 
\bar L(x,y)= \big( \mu x_1,\,\mmt\,\ol x,\, \wmu y_1,\, \mmt\, \ol y\big)+  \fc_\cB. 
\]
Changing  variables $(x,y)\to \bar L(x,y)$ gives
\Be 
\label{op-norm3}
\|\mathcal O_\lambda [\Phi_s, \tilde A_s]\|_{2\to 2}  =(\lambda 2^{-3l})^{-1} (\mu\wmu)^{\frac d2} 
\big\|\mathcal  O_{\smu\wmu\lambda} [\, \bar\Phi,  \bar B\,]\big\|_{2\to 2}. 
\Ee
 
 Freezing  $z_1:=(x_1, y_1)$, we define 
\[
\ol{\mathcal  O}_{z_1}^{\lambda} h(\ol x) = \int e^{i\lambda\bar\Phi_{z_1}(\ol x,\ol y)} \bar B_{z_1} (\ol x, \ol y)h(\ol y)\, d\ol y,
\]
where $\bar \Phi_{z_1}(\ol x,\ol y)=\bar\Phi (x,y)$ and  $\bar B_{z_1} (\ol x, \ol y)=\bar B  (x, y)$. 
Note  $2^{-l}\ge \epsilon_1\wmu^{1/2}$ and $\lambda\wmu^{3/2}\gtrsim 1$.
Since $\mathcal  O_{\lambda} [\, \bar\Phi,  \bar B\,] f=\int   \ol{\mathcal  O}_{(x_1, y_1)}^{\lambda}  f(y_1, \cdot) dy_1$, by  \eqref{op-norm3}  the  estimate \eqref{key-l22} follows if we show
\Be
\label{l2d-1}
\big\| \mathcal O_{z_1}^{\lambda} \|_{2\to 2}\le C \lambda^{-\frac {d-1}2}, \quad  \lambda>0
\Ee
for a constant $C$.  Note $\supp \bar B_{z_1}\subset \mathbb B_{d-1}(0,1)\times \mathbb B_{d-1}(0,1)$.
Moreover,  by  \eqref{bdbs} and \eqref{matdet2} we have $\partial_{\ol x}^\alpha \bar{B}_{z_1} = O(1)$ and $ \det \partial_{\ol x}  \partial_{\ol y}^\intercal   \bar \Phi_{z_1} \sim 1$. Lemma 
\ref{lem:ph-derivative} and \ref{lembd},  as before (cf. \emph{Proof of Lemma \ref{key-osc1}}), give
\[  \partial_{\ol x}^\alpha \partial_{\ol y}^\beta\bar \Phi_{z_1} = O(1), \quad |\alpha|\ge 1,\,  |\beta|=1, 2 \]
 whenever $(\ol x,\ol y)\in \mathbb B_{d-1}(0,1)\times \mathbb B_{d-1}(0,1)$.   Therefore, applying  Lemma \ref{generalized-small2} (with $\omega=1$ and $d$ replaced by $d-1$) to $\mathcal  O_{z_1}^{\lambda} $, we  obtain  
$\eqref{l2d-1}$. 
\end{proof}

\begin{proof}[Proof of Lemma \ref{bdbs}]  We first consider  the case $\bar\alpha=0$. 
Note that $\wt{\psi}(s)A_s(x,y)$ and $\partial_s(\wt{\psi}(s)A_s(x,y))$ are uniformly bounded  on  $\cB\times \tilde\cB$. 
Recall
$|\partial_s^2\Phi_s(x,y)|\lesssim 2^{-3l}$ (see \emph{Proof of Proposition \ref{mm2}}). 
Thus, by \eqref{lower-} we get \eqref{e_bdBs} for $\bar\alpha=0$.

Let  $\bar\alpha\neq 0.$   We have \eqref{est-amp},  and   $\partial_{x}^\alpha \partial_s A_s(x,y)=O((\mu\wmu)^{-{|\alpha|}/{2}})$, which follows by \eqref{est-partial} since 
$2^{-l}\lesssim \mu^{1/2}$.  Therefore,  for \eqref{e_bdBs}
 we need only  to show
\begin{align}
\label{Phis2}
|\partial_{x}^\alpha \partial_s^2\Phi_s(x,y)| & \lesssim 2^{-3l} \mu^{- {|\alpha|}},  
\\
\label{Phis1}
|\partial_{\ol x}^{\ol \alpha} \partial_s \Phi_s(x,y)| & \lesssim 2^{-3l}(\mu\wmu)^{-\frac{|\ol \alpha|}{2}}.  
\end{align}
 
 We verify \eqref{Phis2}  first. Since $\partial_s^2\Phi_s=2^{-2l} \partial_s^2 \mathcal P(x,y, \sxyl)$, it suffices to show 
\Be\label{strong}
|\partial_x^\alpha (\partial_s^2\mathcal P(x,y,S_c^l))|\lesssim 2^{-l}\mu^{-|\alpha|},  \quad \alpha\neq 0.
\Ee
Using \eqref{ph-dd2},   we  write
\[
\frac{\partial^2_s\mathcal P(x,y,S_c^l) }{\inp xy}= \frac{(\cos S_c^l-\cos S_c)^2}{\sin^3 S_c^l}+ (\cos S_c^l-\cos S_c)\frac{(\cos S_c - \tau^+)}{\sin^3 S_c^l}.
\]
Note $   \partial_{x}^\alpha (\sin S_c^l)^{-3}=O(\mu^{-\frac32}\mu^{- {|\alpha|}})$ (see  Remark \ref{derivs}). By the Leibniz rule, \eqref{strong} follows if we show 
\begin{align}
| \partial_{ x}^\alpha(\cos S_c-\tau^+)|&\lesssim  \mu^{1-|\alpha|},
\label{costau}
\\
\label{diffsc}|\partial_x^\alpha(\cos S_c^l - \cos S_c)|&\lesssim 2^{-l}\mu^{\frac12-|\alpha|}. 
\end{align}

 The estimate \eqref{costau}  is clear from \eqref{cossin} since  $\partial_{ x}^\alpha(|x-y|)= O(\mu^{1-|\alpha|})$ for  $(x,y)\in \cB\times \tilde\cB$. 
  To show \eqref{diffsc},  we observe $\partial_{ x}^\alpha(\cos S_c^l-\cos S_c)$ is  given by a linear combination of  the terms
\[
\textstyle (\sin S_c^l-\sin S_c) \prod_{m=1}^{2l-1}\partial_{x}^{a_m} S_c, \qquad (\cos S_c^l-\cos S_c) \prod_{m=1}^{2l} \partial_{x}^{b_m} S_c,  \quad l\ge 0,
\]
where $a_1, \dots, a_{2l-1}, b_1,\dots, b_{2l}\neq 0$;  $|a_1|+\cdots+|a_{2l-1}| =|b_1|+\cdots+|b_{2l}| = |\alpha|$.
Since $\sin S_c^l-\sin S_c = O(2^{-l})$ and $\cos S_c^l-\cos S_c = O(2^{-l}\mu^{1/2})$,  \eqref{diffsc} follows by \eqref{est-partial}.

We now show \eqref{Phis1}. Recalling $\partial_s\Phi_s=2^{-l} \partial_s \mathcal P(x,y, \sxyl)$, 
we have
\[
\partial^{\bar\beta}_{\ol x}\partial_s\Phi_s = 2^{-l}\partial_{\ol x}^{\bar\beta} \partial_s\mathcal P(x,y,S_c^l) + 2^l\partial^2_s\Phi_s\partial^{\bar\beta}_{\ol x} S_c,\quad |\bar\beta| = 1.
\]
Let $\bar\alpha'+\bar\beta=\bar\alpha.$
By \eqref{Phis2}  and \eqref{est-partial}, we have
$\partial_{\bar x}^{\ol \alpha'}(2^l\partial^2_s\Phi_s\partial^{\bar\beta}_{\ol x} S_c)   \lesssim   2^{-3l}(\mu\wmu)^{-{|\ol \alpha|}/{2}}$  since $2^{-l}\ge \varepsilon_1 \wmu^{1/2}$. 
To handle the first term, from \eqref{ph-d} we note
\Be 
\label{j-partial}
\partial_{x_j}\partial_s\mathcal P(x,y,S_c^l)  = \frac{y_j\cos S_c^l - x_j}{\sin^2 S_c^l}. 
\Ee
Thus, it suffices to show  $|\partial_{\ol x}^{\bar\beta}(y_j\cos\sxyl - x_j)|\lesssim (\mu\wmu)^{{(1-|\bar\beta|)}/{2}}$, $2\le j\le d$. 
Since $|y_j|, |x_j| \lesssim (\mu\wmu)^{1/2}$, $2\le j\le d$, using  \eqref{est-partial} one can easily show  the desired bounds  (see Remark \ref{derivs}). 
\end{proof}

The rest of this section is dedicated to proving Lemma \ref{matlem}. 

\subsection{Proof of Lemma \ref{matlem}} 
In order to prove \eqref{erbd},
we start by removing an insignificant  part of $\mathfrak M$.
 Differentiating   \eqref{Ps}, we have 
\begin{equation}
\label{partialxy}
\begin{aligned}
&\quad \partial_x\partial_y^{\intercal}   \Phi_s(x,y) =  \partial_x \partial_y^{\intercal} \mathcal P(x,y, \sxyl)
+
\partial_x S_c \partial_y^{\intercal}  \partial_s \mathcal P(x,y,  \sxyl) +
\\[2pt]
 \partial_x&\partial_s\mathcal P(x,y,  \sxyl) \partial_y^{\intercal}  S_c 
 +\partial_s^2\mathcal P(x,y,  \sxyl) 
\partial_x S_c  \partial_y^{\intercal}  S_c +  \partial_s\mathcal P(x,y,  \sxyl) 
 \partial_x\partial_y^{\intercal}  S_c.\!\!\!\!\!\!\!
\end{aligned}
\end{equation} 
The last term on the right hand side is negligible.

\begin{lem}\label{error1}
Let $s\in \mathbb J$ and $(x,y) \in \cB\times {\tilde \cB}$. Suppose that   \eqref{mucon} and \eqref{interval} hold.  Then, for $\alpha\in \N_0^{d}$
\[
\sin S_c^l\, \partial_x^\alpha \big(\partial_s\mathcal P(x,y,\sxyl) \partial_x\partial_y^{\intercal}  S_c\big) = O(E_1^\alpha).
\]
\end{lem}

To prove Lemma \ref{error1}, we  make use of the following identities: 
\Be
\label{xyS}
 \partial_x S_c(x,y)=\bc(2x-\ba y), \qquad \partial_y^\intercal S_c(x,y)=\bc(2y^\intercal-\ba x^\intercal),
 \Ee
where 
\Be
\label{def-a-c}
\quad  
 \ba(x,y)= \frac{|x|^2+|y|^2}{\inp xy}, 
     \qquad   \bc(x,y)=\frac{\cu(x,y)}{\su(x,y) |x+y||x-y|}. 
\Ee

\begin{proof}[Proof of  \eqref{xyS}]
Differentiating  \eqref{sc}, followed by a simple computation, yields  
\[
\sin S_c(x,y)  \partial_x  S_c(x,y)= \mathrm F(x,y)      x -   2^{-1}  \ba(x,y)   \mathrm F(x,y)  y,
\]
where 
\begin{align*}
 \mathrm F(x,y)=  \  \frac{|x|^2+|y|^2-|x+y||x-y|}{\inp xy |x+y||x-y|}. 
\end{align*}
Here,  we also use the identity  $|x+y|^2|x-y|^2=(|x|^2+|y|^2)^2-4\inp xy^2.$ 
 By  \eqref{sc-} and  \eqref{roots}  we see
$ \mathrm F(x,y)= {2\cos S_c(x,y)}/{|x-y||x+y|}. $
  We thus get  the first identity in \eqref{xyS}.  Since $S_c(y,x)=S_c(x,y)$, the second one  follows from the first by interchanging the roles of $x,y$.
\end{proof}

\begin{proof}[Proof of Lemma \ref{error1}]
Since $|\cos S_c^l - \cos S_c|\sim 2^{-l} \mu^{1/2}$,   \eqref{haha}
gives $|\inp xy-\cos \sxyl(x,y)|\lesssim  \vepc \mmt +2^{-l}\mu^{1/2}.$ 
Thus,  by \eqref{ph-d} and \eqref{small-det}  it follows that  
\[
  \partial_s\mathcal P(x,y,  \sxyls) 
  =O( \vepc  {\tilde\mu}+2^{-2l}) .    
\]
 Combining this and  \eqref{est-partial},   we  
get the desired estimate for $\alpha= 0$. 

We now assume $|\alpha|\ge 1$. 
Thanks to \eqref{est-partial}, it suffices  to show 
\Be 
\label{alpha-}
\partial_x^\alpha\big(\partial_s\mathcal P(x,y,S_c^l)\big)=O\big((\wmu +2^{-2l})\mu^{-\alpha_1} (\mu\wmu)^{-\frac{|\ol\alpha|}2}\big).\Ee
Since $\partial_{x_j}(\partial_s\mathcal{P}(x,y,S_c^l))= \partial_{x_j}\partial_s\mathcal{P}(x,y,S_c^l) + \partial_s^2\mathcal{P}(x,y,S_c^l) \partial_{x_j} S_c$, using \eqref{j-partial}, \eqref{ph-dd2}  and the first  identity in \eqref{xyS},  we write 
\[
    \partial_{x_j}(\partial_s\mathcal{P}(x,y,S_c^l)) =   \frac{  \mathcal I_j^1}{\sin^2 S_c^l} + \frac{\cos S_c^l - \cos S_c}{\sin^3 S_c^l} \,\CI_j^2,
\]
where
\begin{align*}
    \CI_j^1 &:= {y_j\cos S_c - x_j}, 
    \\[2pt]
    \qquad 
    \CI_j^2 &:= {y_j\sin S_c^l - \bc\inp xy(\cos S_c^l - \tau^+)(2 x_j-\ba y_j)}.
\end{align*}

We may assume $\alpha_j\neq 0$ for some $j$. Thus, 
\[ \partial_x^\alpha\big(\partial_s\mathcal P(x,y,S_c^l)\big)=
\partial_x^{\beta}  \Big(  \frac{ \mathcal I_j^1}{\sin^2 S_c^l} \Big )+ \partial_x^{\beta} \Big(  \frac{\cos S_c^l - \cos S_c}{\sin^3 S_c^l} \,\CI_j^2\Big),\] 
where $\beta=\alpha-e_j$. 
By the Leibniz rule, \eqref{diffsc}, and  the bounds in Remark \ref{derivs},    the desired estimate  \eqref{alpha-} follows if we  show  
 \begin{align}
\label{11-1}
    &|\partial_x^\beta \mathcal I_j^1|\lesssim 
             \begin{cases}  
             \wmu\mu^{-\beta_1}(\mu\wmu)^{-{|\ol\beta|}/{2}},  & \qquad   \quad  j=1, 
             \\
             (\mu\wmu)^{1/2}\mu^{-\beta_1}(\mu\wmu)^{-{|\ol\beta|}/{2}},  &   \ \  \qquad \,\, j\neq 1,
             \end{cases}         
    \\
    \label{11-2} 
    &|\partial_x^\beta \mathcal I_j^2 |\lesssim 
    \begin{cases}  
             (\wmu^{1/2} +2^{-l}) \mu^{-\beta_1}(\mu\wmu)^{-{|\ol\beta|}/{2}},  & \quad  j=1, 
             \\
             \mu^{1/2}\mu^{-\beta_1}(\mu\wmu)^{-{|\ol\beta|}/{2}},  &\quad  j\neq 1. 
             \end{cases}  
\end{align}
It should be noted that  slightly weaker bounds are good enough for our purpose  when $\alpha_j\neq 0$ for $j\neq 1$, i.e., $|\ol \alpha|-1=|\ol \beta|$,  since 
there is an improvement of factor $(\wmu/\mu)^{1/2}$ thanks to the particular form of the estimate \eqref{alpha-}.

One can easily  show  \eqref{11-1} for  $j\neq 1$  using \eqref{est-partial}  since $|x_j|, |y_j|\lesssim   (\mu\wmu)^{1/2}$. 
To show \eqref{11-1} for $j=1$, we write  
$
\mathcal I_1^1 
=x_1^{-1}( \inp xy \cos S_c-|x|^2+ \inp {\ol x}{\ol x- {\ol y}\cos S_c})
$. By \eqref{sc-}, we have \begin{align*}
\mathcal I_1^1=\frac{1}{x_1}\Big(\frac{2( \inp xy^2-|x|^2|y|^2)}{|y|^2-|x|^2+|x+y||x-y|} +\inp {\ol x}{\ol x- {\ol y}\cos S_c} \Big).
\end{align*}
Using \eqref{11-1}  with  $j\neq 1$, we have  $\partial_x^\beta(\inp {\ol x}{\ol x- {\ol y}\cos S_c}) = O((\mu\wmu) \mu^{-\beta_1}(\mu\wmu)^{-{|\ol\beta|}/{2}})$.
One can easily   check $\partial_x^\beta(|x|^2|y|^2 - \inp xy^2) = O(\mu\wmu\mu^{-\beta_1}(\mu\wmu)^{-{|\ol\beta|}/{2}})$. 
Also, note $\partial_x^\beta(|y|^2-|x|^2+|x+y||x-y|) = O(\mu^{1-|\beta|})$. 
Since $ |y|^2-|x|^2+|x+y||x-y|\sim \mu$, combining those estimates with the Leibniz rule, we get \eqref{11-1} with $j=1$.

For the proof of  \eqref{11-2}, we first  claim  that we may replace  ${\mathcal I}_j^2$ with  
\[ 
 \tilde {\mathcal I}_j^2 =y_j\sin S_c - \bc\inp xy(\cos S_c - \tau^+)(2 x_j-\ba y_j).\]
 In what follows, we frequently use this type of argument to replace less favorable terms by allowing acceptable errors.
To show the claim, we make use of some easy estimates.  We first note  
\Be
\label{diffsc2}
|\partial_x^\beta(\sin S_c^l - \sin S_c)|\lesssim 2^{-l}\mu^{-|\beta|},  \quad \beta\in \mathbb N^d,
\Ee
which one can show in the same manner as  \eqref{diffsc}. Using  \eqref{est-partial} and the bounds in Remark \ref{derivs}, we also have  
\Be\label{bd-c}
\partial_x^\beta \bc = O(\mu^{-\frac32-|\beta|}),   \quad \beta\in \mathbb N^d. 
\Ee
Noting  $\ba = 2 + |x-y|^2/\inp xy $, via a routine computation we  have
\Be
\label{baba}
    |\partial_x^\beta(2x_i - \ba y_i)|, \ |\partial_x^\beta(2y_i - \ba x_i)|
    \lesssim\begin{cases}
    \mu\mu^{-\beta_1}(\mu\wmu)^{-{|\ol\beta|}/{2}}, & i = 1, 
    \\[2pt]
    (\mu\wmu)^{1/2}\mu^{-\beta_1}(\mu\wmu)^{-{|\ol\beta|}/{2}}, & i\ge 2,
    \end{cases}
\Ee
for $\beta\in \mathbb N^d$.  Indeed, to see this  one needs only to write  $2x_i - \ba y_i=(2\inp xy x_i - (|x|^2+|y|^2)\ba y_i)/\inp xy$, and 
$2y_i - \ba x_i$ can be handled similarly.  Putting together the estimates \eqref{diffsc2},  \eqref{diffsc},  \eqref{costau}, \eqref{bd-c}, and \eqref{baba}, we now  see 
\[
\partial_x^\beta  ( \mathcal I_j^2-\tilde {\mathcal I}_j^2)=O\big( 2^{-l}  \mu^{-\beta_1}(\mu\wmu)^{-{|\ol\beta|}/{2}}\big),\quad \beta\in \mathbb N^d.
\]
In view of \eqref{11-2}, the difference is an acceptable error. Therefore, as claimed above,  we only have to show 
\eqref{11-2} for   $\tilde{\mathcal I}_j^2$ replacing   ${\mathcal I}_j^2$. 

To do this, using the first equality in   \eqref{cossin} and \eqref{def-a-c},  we note
\[
\sin S_c\, \tilde{\mathcal I}_j^2 = \bar {\mathcal I}_j :=y_j\sin^2 S_c + \cos S_c(2x_j-\ba y_j). 
\]
Thus, the desired bound  \eqref{11-2} follows once we have
\Be
\label{11-22} 
    |\partial_x^\beta \bar{\mathcal I}_j |\lesssim 
    \begin{cases}  
             (\wmu^{1/2} +2^{-l}) \mu^{1/2-\beta_1}(\mu\wmu)^{-{|\ol\beta|}/{2}},  & \quad  j=1, 
             \\
             \mu \mu^{-\beta_1}(\mu\wmu)^{-{|\ol\beta|}/{2}},  &\quad  j\neq 1. 
             \end{cases}       
\Ee

Since  $|\ol x|, |\ol y|=O((\mu\wmu)^{1/2}),$ \eqref{11-22} for  $j\neq 1$ is easy. Indeed, it follows by \eqref{baba} and the bounds in Remark \ref{derivs}. 
To show \eqref{11-22} with $j=1$, we break 
\[  \bar{\mathcal I}_1=  \bar{\mathcal I}_{1,1}- \bar{\mathcal I}_{1,2},  \] 
where 
\begin{align*}
 \bar{\mathcal I}_{1,1}&:= x_1^{-1}(\inp xy \sin^2 S_c + \cos S_c(2|x|^2-\ba \inp xy)),\\
 \bar{\mathcal I}_{1,2}&:=  x_1^{-1}(\inp{\ol x}{\ol y}\sin^2 S_c + \cos S_c(2|\ol x|^2 - \ba \inp{\ol x}{\ol y})).
\end{align*}
Since  $|\ol x|, |\ol y|=O((\mu\wmu)^{1/2}),$ one can easily see $|\partial_x^\beta \bar{\mathcal I}_{1,2}|\lesssim \mu\wmu \mu^{-\beta_1}(\mu\wmu)^{-{|\ol\beta|}/{2}}$.  
To handle $\bar{\mathcal I}_{1,1}$,  using $\ba \inp xy = |x|^2+|y|^2$,  the second equality in \eqref{cossin}, and $|x-y||x+y|=|x|^2+|y|^2-2\inp xy \cos S_c$, successively, we note  \[
\bar{\mathcal I}_{1,1} = 2 x_1^{-1} \cos S_c\inp{x}{x-\cos S_c\, y}.\]  
By \eqref{11-1}, we have $\partial_x^\beta (x-\cos S_c\, y)=O((\mu\wmu)^{1/2}\mu^{-\beta_1}(\mu\wmu)^{-{|\ol\beta|}/{2}})$.  Therefore, we get  \eqref{11-22} for $j=1$ similarly as before (see Remark \ref{derivs}). 
\end{proof}

Let us set
\begin{align*}
\mathbf G(x,y)&=
    -\sin\sxyl\big(\,\partial_x S_c  \partial_y^{\intercal}  \partial_s \mathcal P(x,y,  \sxyl) 
+ \partial_x\partial_s\mathcal P(x,y,  \sxyl) \partial_y^{\intercal}  S_c  \,\big),\nonumber
 \\
\mathbf F(x,y)&= -\sin\sxyl\partial^2_s\mathcal P(x,y,\sxyl)\partial_x S_c\partial_y^{\intercal} S_c,
\\ 
 \mathbf K(x,y)&= \  \mathbf I_d+ \mathbf G(x,y) + \mathbf F(x,y).
\end{align*}
In order to prove Lemma \ref{matlem}, by \eqref{partialxy} and Lemma \ref{error1} it is sufficient to show that the entries $\tilde{\mathcal E}_{i,j}$ of 
the matrix 
\[ \tilde{\mathcal E}:=\mathbf K-\mathfrak M^0 \] satisfy  \eqref{erbd} in place of $\CE_{i,j}$.
A simple  computation gives  
\begin{align*}
    \partial_x\partial_s\mathcal P(x,y, \sxyl)=\frac{\cos  \sxyl\, y-x}{\sin^2  \sxyl }, \qquad 
     \partial_y^\intercal \partial_s\mathcal P(x,y, \sxyl)=\frac{\cos  \sxyl \, x^\intercal-y^\intercal}{\sin^2 \sxyl}.
\end{align*}
Using those identities, we write 
\begin{align}
\label{FG}
    \mathbf G 
    &= -\frac \bc\sul  \Big((2\cul+\ba)(xx^\intercal+yy^\intercal)-4xy^\intercal-2\ba\cul\, yx^\intercal \Big),
    \\[2pt]
    \label{FG2}
    \mathbf F     &=\sul \, \partial_s^2\mathcal P (x,y, \sxyl) \, \bc^2\Big(\,2\ba(xx^\intercal+yy^\intercal)-4 xy^\intercal-\ba^2 yx^\intercal\,\Big).
\end{align}

In what follows, for $i=1,2$,  we denote $E(x,y)=\dO(E_i)$ if  $\partial_x^\alpha E(x,y)=O(E_i^\alpha)$ for any $\alpha$ and $(x,y)\in \cB\times \tilde\cB$.
We first show \eqref{erbd} for $(i,j)=(1,1)$, which is more involved than the others.

\subsubsection*{Proof of \eqref{erbd} for $(i,j)=(1,1)$}  We consider $\mathbf G_{1,1}$ and $\mathbf F_{1,1}$, from which we  discard  some harmless parts. By \eqref{FG},
\[
    \mathbf G_{1,1} 
    = -\frac{\bc}{\sin S_c^l}\Big((\ba + 2)(x_1-y_1)^2 + 2(\cos S_c^l - 1)(x_1^2+y_1^2 - \ba\, x_1 y_1)\Big). 
    \] 
Using $\ba \inp xy = |x|^2+|y|^2$, we  have      
  \[  \mathbf G_{1,1}  = -\frac{\bc}{\sin S_c^l}(\ba + 2)(x_1-y_1)^2 +2\bc\,\frac{ (1-\cos S_c^l) }{\sin S_c^l} (\ba\inp{\ol x}{\ol y}-|\ol x|^2-|\ol y|^2).
\] 
The second term, which we denote by $ \tilde {\mathbf G}_{1,1}$,  is $\dO(E_1)$. Indeed, since  $\partial_x^\alpha (\ba \inp{\ol x}{\ol y}-|\ol x|^2-|\ol y|^2)\big) =O(\mu\wmu\mu^{-\alpha_1}(\mu\wmu)^{-{|\ol\alpha|}/{2}})$,   the bounds in Remark \ref{derivs} and  \eqref{bd-c} show $ \partial_x^\alpha \tilde {\mathbf G}_{1,1}=O(\wmu \mu^{-\alpha_1}(\mu\wmu)^{-{|\ol\alpha|}/{2}})$. 
Thus, by the first identity in \eqref{def-a-c} we get
\[
\mathbf G_{1,1} =-\frac{\cos{S_c}(x_1-y_1)^2}{\sin{S_c}\sin{\sxyl}|x-y|}\frac{|x+y|}{\inp xy} + \dO( E_1). 
\]
By writing $(x_1-y_1)^2=|x-y|^2-|\ol x-\ol y|^2$ and using  \eqref{cossin} (the second equality),   the first term on the right hand side equals 
\[     \frac{|\ol x-\ol y|^2}{|x-y|^2} + |\ol x-\ol y|^2\frac{\sin S_c - \sin S_c^l}{|x-y|^2{\sin S_c^l}}-\frac{|x+y||x-y| \cos \sxy}{\inp xy \sul\su}. \]  
Denote the second term by $\bar {\mathbf G}_{1,1} (x,y)$. 
By \eqref{diffsc2} and \eqref{est-partial}, we have $\partial_x^\alpha\bar {\mathbf G}_{1,1}=O(2^{-l}\wmu \mu^{-3/2}\mu^{-\alpha_1}(\mu
\wmu)^{-|\bar\alpha|/2})$, so $\bar {\mathbf G}_{1,1}=\dO(E_1)$ by \eqref{interval}. Thus, we obtain 
\[
\mathbf G_{1,1} -\mathfrak M^0_{1,1}  + \dO(E_1)
= -\frac{|x+y||x-y| \cos \sxy}{\inp xy \sul\su} = -\frac{(\tau^+-\cu) \cos \sxy}{ \sul\su} .
\]
We use \eqref{cossin} for the second equality.

From \eqref{FG2}, we have 
\begin{align*}
    \mathbf F_{1,1} &= \sin S_c^l \partial_s^2\mathcal P(x,y,S_c^l)\,\bc^2 \big(2\ba(|x|^2+|y|^2)-(4+\ba^2)\inp xy) \\
    &\qquad -\sin S_c^l \partial_s^2\mathcal P(x,y,S_c^l)\,\bc^2(2\ba(|\ol x|^2+|\ol y|^2)-(4+\ba^2)\inp{\ol x}{\ol y}).
\end{align*}
We denote the second term by $\tilde {\mathbf F}_{1,1}$.  By \eqref{strong} and \eqref{bd-c}, we see 
$\partial_x^\alpha \tilde {\mathbf F}_{1,1}=O(2^{-l}\wmu \mu^{-3/2}\mu^{-\alpha_1}(\mu\wmu)^{-{|\ol\alpha|}/{2}})$, 
so  $ \tilde {\mathbf F}_{1,1}=\dO(E_1)$. Meanwhile, \eqref{def-a-c}  gives
\[
2\ba(|x|^2+|y|^2)-(4+\ba^2)\inp xy = \frac{|x-y|^2|x+y|^2}{\inp xy} = \bc^{-2}\frac{\cos^2 S_c}{\inp xy \sin^2 S_c}.
\]
Combining this with \eqref{ph-dd2} yields
\begin{align*}
    \mathbf F_{1,1} 
    =\frac{ (\tau^+-\cul) (\cul-\cu) \cos^2\sxy}{\sult\sut}+\dO(E_1).
\end{align*}

We now set 
\begin{align*}
D_1&=   (\tau^+-\cul)(\cul-\cu)\cos^2 S_c, 
\\
D_2&= (\tau^+-\cu) (\sul-\su)\sul\cu, 
\end{align*}
and 
\[ 
A= \frac{D_1+  D_2 }  {\sut\sult}. \]
Since $ \mathbf K_{1,1}=1+\mathbf F_{1,1}+\mathbf G_{1,1}$, using the above equalities,  we obtain  
    \begin{align*}
    \mathbf K_{1,1}-\mathfrak M^0_{1,1}& = A +\dO(E_1). 
       \end{align*}
Indeed,  from \eqref{sc-} and \eqref{roots}  note $\sut=(\tau^+-\cu)\cu$, which then  gives   $D_2=\sut\sult-{(\tau^+-\cu)}  \su\sul   \cos \sxy.$ 

To complete the proof, it remains to show  $A=\dO(E_1).$  Let  $\tilde D_1$ denote  $D_1$ in which 
the first $\cul$ is replaced with $\cu$,  and  by $\tilde D_2$ we denote  $D_2$ in which the second $\sul$ is replaced with $\su$.  
Then,   using \eqref{diffsc} as before, we see $\partial_x^\alpha( D_1- \tilde D_1)=O(2^{-2l}\mu\mu^{-\alpha_1}(\mu\tmu)^{-|\ol \alpha|/2})$.
Similarly, $\partial_x^\alpha( D_2- \tilde D_2)=O(2^{-2l}\mu\mu^{-\alpha_1}(\mu\tmu)^{-|\ol \alpha|/2})$ by \eqref{diffsc2} and \eqref{costau}. 
So, $( D_1-\tilde D_1+  D_2-\tilde D_2)/(\sut\sult)=\dO(E_1)$.  Thus, we have 
\[ A= \frac{\tilde D_1+  \tilde D_2 }  {\sut\sult}   +\dO(E_1).\]

Elementary trigonometric identities
give $\tilde D_1+\tilde D_2= 2\cu(\cu-\tau^+) \sin^2{(({\sxyl-S_c})/{2})}$. 
Therefore,    
\begin{align*}
A =          
 2  \frac{\cu(\cu-\tau^+)}{\sut\sult}  \sin^2{\Big(\frac{\sxyl-S_c}{2}\Big)} + \dO(E_1).
\end{align*}
Since $S_c^l-S_c = 2^{-l} s$, using  \eqref{costau} and the bounds in Remark \ref{derivs}, we conclude $A=\dO(E_1)$.
\qed

\medskip

Before we begin to  prove  \eqref{erbd} for $(i,j)\neq (1,1)$, we show that the contribution of  $\mathbf F$ is negligible. By   \eqref{FG2} we have 
\begin{align*}
    \mathbf F_{i,j} 
        = -\sin S_c^l\partial_s^2\mathcal P(x,y,S_c^l)\bc^2(2 x_i-\ba y_i)(2y_j - \ba x_j).
\end{align*}
If $i,j\ge 2$,  using  \eqref{strong},  \eqref{baba}, and \eqref{bd-c} (also see Remark \ref{derivs}), we  get
$ \partial_x^\alpha \mathbf F_{i,j}=  O(2^{-l}\wmu \mu^{-3/2}\mu^{-\alpha_1}(\mu
\wmu)^{-|\ol \alpha|/2})$, which shows  $\mathbf F_{i,j}=\dO( E_1)$. If  $i = 1, j\ge 2,$ or $j=1 , i\ge 2$, we  obtain  
$ \partial_x^\alpha \mathbf F_{i,j}=  O(2^{-l}\wmu^{1/2} \mu^{-1}\mu^{-\alpha_1}(\mu
\wmu)^{-|\ol\alpha|/2})$ similarly, so it follows that  $\mathbf F_{i,j}=\dO( E_2)$. 

Therefore, the following completes  the proof of Lemma 4.6:
 \begin{equation}\label{Gij}
\mathbf G_{i,j} -\mathfrak M_{i,j} ^0= 
\begin{cases}  \dO( E_1),   \qquad\qquad \quad \quad  i,j\ge 2,
       \\[2pt]
          \dO( E_2),  \qquad \quad  i = 1, j\ge 2,\text{ or } j=1 , i\ge 2.
           \end{cases}
       \end{equation}

\subsubsection*{Proof of  \eqref{Gij}} We consider the case $i,j\ge 2$ first. By \eqref{FG}
we have
\[
\mathbf G_{i,j} = - \frac{\bc}{\sin S_c^l}\big((2\cos S_c^l + \ba)(x_ix_j+y_iy_j) - 4x_iy_j -2 \ba \cos S_c^l\, y_ix_j\big).
\]
By ${\mathbf G}_{i,j}'$ we denote $\mathbf G_{i,j}$ in which  we replace $\mathfrak A$ with $2$. 
Since $\ba = 2 + |x-y|^2/ \inp xy $, $\partial_x^\alpha (\ba-2)=O(\mu^{2-|\alpha|})$. Using \eqref{bd-c}, 
we see  $\partial_x^\alpha(\mathbf G_{i,j}-{\mathbf G}_{i,j}')=O(\wmu\mu^{1-\alpha_1}(\mu
\wmu)^{-|\ol\alpha|/2})$ 
because $ x_i, x_j, y_i, y_j=O((\mu\wmu)^{1/2})$, $i,j\ge 2$. 
Thus, it is enough to consider ${\mathbf G}_{i,j}'$. 
Furthermore, by \eqref{diffsc} and  \eqref{bd-c},  $\cos \sxyl$ in ${\mathbf G}'_{i,j}$  can similarly be replaced with $\cos S_c$ if we allow an error of $\dO(E_1)$. Thus, 
\Be
\label{GG}
\mathbf G_{i,j}
= \tilde{\mathbf G}_{i,j}  + \bar{\mathbf G}_{i,j}+\dO(E_1),  \quad   i,j\ge 2, 
\Ee
 where 
 \vspace{-3pt}
\begin{align*}
  \tilde{\mathbf G}_{i,j}=  &-\frac{2\bc}{\sin S_c^l}(x_j-y_j)\big(2x_i-(1+\cos S_c)y_i\big), 
    \\
\bar{\mathbf G}_{i,j}=    &\frac{2\bc}{\sin S_c^l}(x_i-y_i)(1-\cos S_c)x_j. 
\end{align*}

It is easy to see $\partial_x^\alpha\bar{\mathbf G}_{i,j}= O(\wmu  \mu^{-\alpha_1}(\mu
\wmu)^{-|\ol\alpha|/2})$  for $i,j\ge 2$, thus  $\bar{\mathbf G}_{i,j}=\dO(E_1)$. Using \eqref{diffsc2}, by $\sin S_c$ we may replace  $\sin S_c^l$ in the expression of $\tilde{\mathbf G}_{i,j}$  with an error of $\dO(E_1)$. 
Then, applying  \eqref{def-a-c} and \eqref{cossin}, we get 
\Be\label{ggg}
 {\mathbf G}_{i,j} =-\frac{2\inp xy}{|x+y|^2} \frac{(x_j-y_j)}{|x-y|^2}\big (2x_i-(1+\cos S_c)y_i\big) + \dO(E_1).
\Ee
As before, we may replace  ${2\inp xy }/{|x+y|^2}$ by $1/2$, since   $\partial_x^\alpha (1/2- 2\inp xy/|x+y|^2)=O(\mu^2 \mu^{-|\alpha|})$, 
$ \partial_x^\alpha (2x_i-(1+\cos S_c)y_i)=O((\mu\wmu)^{1/2} \mu^{-\alpha_1}(\mu
\wmu)^{-|\ol\alpha|/2})$
and $\partial_x^\alpha (x_j-y_j)=O((\mu\wmu)^{1/2} \mu^{-\alpha_1}(\mu
\wmu)^{-|\ol\alpha|/2})$, $i, j\ge 2$. 
Therefore,   we have 
\[ 
 {\mathbf G}_{i,j} = \mathfrak M_{i,j}^0 +\frac{(x_j-y_j)(\cos S_c-1)y_i }{2|x-y|^2}+ \dO(E_1).
\] 
Note  $ \partial_x^\alpha ((x_j-y_j)(1-\cos S_c)y_i/ |x-y|^2)=O( \wmu \mu^{-\alpha_1}(\mu\tmu)^{-|\ol \alpha|/2})$. Thus,  \eqref{Gij} follows  for $i,j\ge 2$. 

We now show \eqref{Gij}  when $i=1$, $j\ge 2$, or  $j=1$, $i\ge 2$. We only consider the case $i=1$, $j\ge 2$,  
since the other one can be handled in the same manner.  Since $x_j, y_j=O((\mu\wmu)^{1/2})$ for $j\ge 2$, repeating the same argument used to show \eqref{GG}, we obtain   
$\mathbf G_{1,j}
= \tilde{\mathbf G}_{1,j}  + \bar{\mathbf G}_{1,j}+\dO(E_2). 
$
Here we use the same notations as above. 
Note $ \partial_x^\alpha (2x_1-(1+\cos S_c)y_1)=O(\mu \mu^{-\alpha_1}(\mu\tmu)^{-|\ol \alpha|/2})$. By this and \eqref{bd-c}, $\partial_x^\alpha\bar{\mathbf G}_{1,j}=O(\mu^{1/2}\wmu^{1/2} \mu^{-\alpha_1}(\mu\tmu)^{-|\ol \alpha|/2}))$, thus ${\mathbf G}_{1,j}=\tilde {\mathbf G}_{1,j}+  \dO(E_2)$. Allowing an error of  $\dO(E_2)$, as before, we may replace $\sin S_c^l$ with $\sin S_c$ in $\tilde {\mathbf G}_{1,j}$. Consequently, using   \eqref{def-a-c} and \eqref{cossin} (cf. \eqref{ggg}), we obtain  
\[
 \tilde {\mathbf G}_{1,j}=-\frac{2\inp xy}{|x+y|^2} \frac{(x_j-y_j)}{|x-y|^2}\big(2x_1-(1+\cos S_c)y_1\big) + \dO(E_2).
\]
We may also replace, as above,  the factor ${2\inp xy }/{|x+y|^2}$ by $1/2$ with an error of $\dO(E_2)$.
Therefore, \eqref{Gij} follows for $i=1$, $j\ge 2$.
\qed

\bigskip

\noindent{\bf Acknowledgements.}  This work was supported by the POSCO Science Fellowship and   NRF-2020R1F1A1A01048520 (E.J.)  and NRF-2021R1A2B5B02-001786 (S.L. $\&$ J.R.). 
The authors thank the referee for careful reading and valuable comments. 

\end{document}